\documentclass{amsart}

      \theoremstyle{plain}
      \newtheorem{theorem}{Theorem}[section]
      \newtheorem{lemma}[theorem]{Lemma}
      \newtheorem{corollary}[theorem]{Corollary}

      \theoremstyle{definition}
      \newtheorem{definition}[theorem]{Definition}

      \theoremstyle{remark}
      \newtheorem{remark}[theorem]{Remark}

      \numberwithin{equation}{section}

      \newcommand{\isom}{\cong}
      \newcommand{\id}{\text{Id}}
      \newcommand{\im}{\text{Im}}
      \newcommand{\abs}[1]{\lvert#1\rvert}
      \newcommand{\innprod}[2]{\langle#1,#2\rangle}
      \newcommand{\norm}[1]{\lVert#1\rVert}

      \newcommand{\bj} {\bar{j}}
      \newcommand{\bk} {\bar{k}}
      \newcommand{\nn} \nonumber

      \newcommand{\bpizb}[1] {\bar{\pi}_{\lambda} (Z_{\overline{#1} \lambda})}
      \newcommand{\bpiz}[1] {\bar{\pi}_{\lambda} (Z_{#1 \lambda})}

\begin{document}

\bibliographystyle{amsplain}

\title{The Laplacian on $p$-forms on the Heisenberg group}

\author{Luke Schubert}
\address{The University of Adelaide, Adelaide, Australia}
\email{lschuber@maths.adelaide.edu.au}

\subjclass{Primary 58G25, 58G18; Secondary 43A80}

\keywords{Novikov-Shubin invariants, $L^2$ cohomology, spectral theory}

\begin{abstract}
The Novikov-Shubin invariants for a non-compact Riemannian manifold $M$
can be defined in terms of the large time decay of the heat operator
of the Laplacian on $L^{2}$ $p$-forms, $\triangle_{p}$, on $M$.

For the $(2n+1)$-dimensional Heisenberg group $H^{2n+1}$, the
Laplacian $\triangle_{p}$ can be decomposed into operators
$\triangle_{p,n}(k)$ in unitary
representations $\bar{\beta}_{k}$ which, when restricted to the centre
of $H$, are characters (mapping $\omega$ to $\exp(-ik\omega)$).
The representation space is an anti-Fock space ($\mathcal{F}^{-k}_n$),
of anti-holomorphic functions $F$ on $\mathbb{C}^n$ such that 
$\int_{\mathbb{C}^n} \abs{F(\bar{z})}^2 e^{-1/4 k \abs{z}^2} dz 
< + \infty$.

In this paper, the eigenvalues of $\triangle_{p,n}(k)$ are calculated,
for all $n$ and $p$,
using operators which commute with the Laplacian;
this information determines the $p$th Novikov-Shubin invariant of $H^{2n+1}$.
Further, some eigenvalues of operators connected with nilpotent Lie groups
of Heisenberg type are calculated in the later sections.
\end{abstract}

\maketitle

\section{Introduction}

The Novikov-Shubin invariants for a non-compact Riemannian manifold $M$
can be defined in terms of the large time decay of the heat operator
of the Laplacian on $L^{2}$ $p$-forms, 
on $M$, which we'll denote by $\triangle_p$.

For the $(2n+1)$-dimensional Heisenberg group $H^{2n+1}$, the
Laplacian can be decomposed into operators $\triangle_{p}(k)$ in irreducible
unitary
representations $\bar{\beta}_{k}$ which, when restricted to the centre
of $H$, are characters (mapping $\omega$ to $\exp(-ik\omega)$).
In this paper, the eigenvalues of $\triangle_{p}(k)$ are calculated;
these determine the Novikov-Shubin invariants of $H^{2n+1}$.

Novikov-Shubin invariants are a relatively new set of topological
invariants, usually defined analytically, of certain non-compact
Riemannian manifolds. They were first defined in \cite{ns2, ns1},
but are most comprehensively discussed in \cite{g&s}, from which
the below definition is taken.

They are related to the $L^{2}$ Betti numbers, in the
following way.
We can define a function $\theta_{p}(t)$, depending on
the manifold $M$, which is a positive function of
$t \in \mathbb{R}_{+}$. Then the $p$th $L^{2}$ Betti number
$b_{p}^{(2)}$ is equal to the limit as $t \rightarrow +\infty$
of $\theta_{p}(t)$, while the $p$th Novikov-Shubin number $\alpha_{p}$
measures the (degree of the inverse polynomial)
rate at which this limit is approached.

The theory of $L^{2}$ torsion (see for example \cite{cfm}) is also
closely linked to that of Novikov-Shubin invariants; for example, if
all the Novikov-Shubin invariants of a manifold are positive, then the
$L^{2}$ torsion of that manifold is defined.  Further, Novikov-Shubin
invariants have been placed in a more abstract, categorical setting
and more naturally linked with torsion and $L^2$ cohomology by Farber
in \cite{farber}.

These invariants do not exist for all manifolds, but generalised
invariants to cover the exceptions have been defined. (See, again,
\cite{g&s}, and also \cite{ccmp, vmdirac}, while a combinatorial 
definition is given in \cite{efremov}, and the Novikov-Shubin invariants
of complexes of Hilbert spaces are defined in \cite{g&s, l&l}.)

For $d$ the usual exterior derivative on square-integrable
      $p$-forms, and $d^*$ its adjoint with respect to the Riemannian
      metric, we define the Laplacian on
      $L^2$ $p$-forms to be
      \[ \triangle_p = dd^* + d^* d. \]
      Note that the Laplacian is a self-adjoint, positive, second order
      elliptic differential operator; the Laplacian on functions,
      $\triangle_0$, is the typical such elliptic differential operator.
      The Laplacian on forms is somewhat more complex, but differs from
      $\triangle_0 \otimes \id$ only in first and zeroth order terms.

We first define the heat operator $e^{-t \triangle_{p}}$ for all $t
>0$ using the spectral theorem for self-adjoint operators. Then for
$\Gamma$ a discrete subgroup
      of the isometry group of $M$, such that $M/\Gamma$ is a compact
manifold, we can define
      a certain von Neumann trace ${\rm Tr}_\Gamma$ on
      the $\Gamma$-invariant operators of $B(L^2(M))$,
      and thus the function $\theta_{p}(t):={\rm Tr}_\Gamma (e^{-t
      \triangle_p})$, mentioned above, can be defined for all positive $t$.
      Then $\theta_{p}(t)$ approaches
      $b_{p}^{(2)}$, the $p$th $L^{2}$ Betti number of $M$, for
      large $t$.
      If furthermore $\theta_{p}(t)-b_{p}^{(2)}$
      is of order $t^{-\alpha_p}$
      and $t^{-\alpha_{p}}$ is of order $\theta_{p}(t)-b_{p}^{(2)}$
      as $t \rightarrow \infty$, then
      we say that $\alpha_p$ is the $p$th Novikov-Shubin invariant
(see \cite{g&s}).

 It is known that $\alpha_{p}$ is independent of the choice of 
 $\Gamma$-invariant metric
 on $M$; 
other properties of $\alpha_{p}$ are discussed in the main
 text.

In this paper, we calculate the Novikov-Shubin invariants
of the $(2n+1)$-dimensional Heisenberg group $H^{2n+1}$.
The method chosen is to examine not only the Laplacian on
$L^{2}$ $p$-forms on $H^{2n+1}$, but also this operator in an
irreducible, unitary
representation of $H^{2n+1}$. We study the spectrum of this latter
operator and thereby derive all the Novikov-Shubin
invariants for each Heisenberg group.

Recall that the Heisenberg group of dimension $(2n+1)$ (hereafter denoted
by $H^{2n+1}$, or $H$ if the dimension is clear) is a $2$-step nilpotent
Lie group. (It arises naturally in
Quantum Mechanics; it is also, in some sense, the simplest
non-abelian nilpotent Lie group.)

We choose a left-invariant metric on $H$; then the Laplacian
$\triangle_{p}$ defined with respect to this metric is (left) $H$-invariant.

      Now since $H$ is a Lie group, its tangent bundle is trivial;
      so $\triangle_p$ can be thought of as a matrix, with entries
      which are
      differential operators on $L^2(H)$. But we know from the abstract
      Plancherel theorem that $L^2 (H)$ splits into a
      direct integral of Hilbert spaces:
      \[ L^2(H) = \int_{\mathbb R}^\oplus
      \mathcal{F}^{k}_{n} \otimes \mathcal{F}^{-k}_{n} |k|^n dk \]
      where $k$ corresponds to the Fock-Bargmann representation
      $\beta_k$ with parameter $k$. The Laplacian $\triangle_p$ also
      splits under this
      direct integral, with the corresponding operator in each
      Hilbert space being denoted by $\triangle_p(k)$.

The central result of this paper is Theorem \ref{T:all_evals},
which lists all the eigenvalues of $\triangle_{p,n}(k)$ for 
$p \le n$ and $k>0$ (though not, 
in general, their multiplicity). In particular, the lowest eigenvalue
of $\triangle_{p,n}(k)$, again for $p \le n$ and $k>0$, is 
$k^{2}+(n-p)k$, which has multiplicity $\binom{n}{p}$. 
Further, Theorem \ref{T:all_evals} implies that the spectrum of
$\triangle_{p,n}(k)$ contains the spectrum of $\triangle_{p-1,n-1}(k)$
for all $p \ge 1, n \ge 2, p \le n$.

Using this theorem, we calculate exactly the Novikov-Shubin invariants
of the Heisenberg group in Corollary \ref{C:NSinv}: 
\begin{eqnarray*}
\alpha_p(H^{2n+1}) = \left\{ \begin{array}{ll}
							n+1, & p \ne n, n+1, \\
							\frac{1}{2}(n+1), & p = n, n+1.
							  \end{array}
					 \right.					
\end{eqnarray*}

The results of Varopoulos in \cite{varopoulos} determine 
$\alpha_0(M)$ explicitly for all manifolds $M$, which agree with the above
for the case $p=0$. Further, Corollary \ref{C:NSinv} refines the 
following inequalities for $\alpha_p(H^{2n+1})$ which were proved in 
\cite{lott}:
\begin{eqnarray*}
\alpha_p(H^{2n+1}) \le \left\{ \begin{array}{ll}
							n+1, & p \ne n, n+1, \\
							\frac{1}{2}(n+1), & p = n, n+1,
							  \end{array}
					 \right.					
\end{eqnarray*}
where our definition of $\alpha_p$ differs by a factor of $2$ from Lott's.

A result analogous to, but weaker than, Theorem \ref{T:all_evals} can
be found in
\cite{g&w}, where the Laplacian on a quotient 
$H/\Gamma$ of the Heisenberg group by a discrete, cocompact subgroup 
$\Gamma$ is considered. There, the eigenvalues of the decomposition of this
operator in characters of $H$ are calculated, rather than the eigenvalues
in an infinite-dimensional representation, as in this paper. 

The algebraic methods which we use to simplify the problem of calculating
eigenvalues of the Laplacian on the Heisenberg group have their analogues
for other nilpotent Lie groups. In the final two sections of this paper,
we generalise these methods to the case of Heisenberg-type groups
and obtain some information on the spectrum of the Laplacian. In 
particular, we obtain estimates of the lowest eigenvalue of the Laplacian
on $1$-forms on a family of nilpotent Lie groups with two-dimensional
centre, and thus calculate the first Novikov-Shubin invariant of these
Lie groups.

The results of this paper form an extension of the results of my thesis
\cite{thesis} which was supervised by Alan Carey and Varghese Mathai;
many thanks are due them for all their patience and encouragement.

\section{Novikov-Shubin invariants}

We define Novikov-Shubin invariants as in \cite{g&s}.

Let $M$ be a non-compact oriented Riemannian manifold on which a discrete
infinite group $\Gamma$ acts freely, such that the quotient $X := M/\Gamma$
is a compact manifold.

Let $\mathcal{A}$ be the algebra of all bounded linear operators on
$L^2(M)$ which commute with the action of $\Gamma$: it can be shown
that $\mathcal{A}$ is a von Neumann algebra \cite{atiyah}.

There is a von Neumann trace on $\mathcal{A}$, denoted by ${\rm Tr}_\Gamma$,
first defined by Atiyah in \cite{atiyah}. If an operator $A \in
\mathcal{A}$ is of $\Gamma$-trace class, has smooth kernel $K_A (x,y)$
(which is a distribution on $M \times M$),
%(which is an $L^2$ function and not just a distribution), 
and is positive
and self-adjoint, then
 \[ {\rm Tr}_\Gamma A = \int_{\mathcal{F}} K_A (x,x) d \mu (x) \]
where $\mathcal{F}$ is the fundamental domain for the action of $\Gamma$
on $M$ and $\mu$ is Haar measure on $\mathcal{F}$. %... measure on F.

Let $\triangle_p$ denote the Laplacian on smooth, compactly supported
$p$-forms on $M$. The closure of this operator, $\tilde{\triangle}_p$,
has domain the first (generalized) Sobolev space on $p$-forms, 
which is dense in
the set of $L^2$ $p$-forms. This space is defined as the closure of 
smooth, compactly
supported $p$-forms on $M$ with respect to the norm
\[ \norm{\omega}_1 = \innprod{(\id + \triangle_p) \omega}
{(\id + \triangle_p) \omega}_2 \]
where $\innprod{.}{.}_2$ is the usual $L^2$ inner product on $p$-forms,
here applied to $\omega$ in the sense of distributions
 (see \cite{atiyah} or \cite{dodziuk}).
Hereafter we write $\triangle_p$ instead of $\tilde{\triangle}_p$ and refer to
this operator as acting on $L^2$ $p$-forms.

Using the spectral theorem for self-adjoint operators, we can form the
operator $e^{-t \triangle_p}$ for all positive $t$.
We can then define a function $\theta_p(t)$ for all $t > 0$, by
\[ \theta_p (t) := {\rm Tr}_\Gamma e^{-t \triangle_p}. \]
It was shown in \cite{g&s} 
that $\theta_p(t) \rightarrow \bar{b}_p$,
the $p$th $L^2$ Betti number (first defined in \cite{atiyah}), as
$t \rightarrow \infty$.

If for some constants $C, t_0$ and $\alpha$, we have
\[ C^{-1} t^{- \alpha} \le \theta_p(t) - \bar{b}_p \le C t^{-\alpha} \]
for all $t > t_0$, we say that $\alpha = \alpha_p(M, \Gamma)$ is
the $p$th Novikov-Shubin invariant of $(M, \Gamma)$.

This was not the original definition of $\alpha_{p}(M)$, but it was
proved to be equivalent in \cite{g&s}. It is the most useful definition
for our purposes.

It has been shown that $\alpha_p(M)$ is invariant of choice of
$\Gamma$-invariant metric, and furthermore is a homotopy
invariant. This last statement was proved in \cite{g&s}, but the proof
is complicated, and relies on assumptions that certain operators are
bounded. For an alternative proof for closed manifolds, which uses
standard topological techniques, see \cite{bmw}.

      \section{The Plancherel theorem for the Heisenberg group}

The Heisenberg group of dimension $2n+1$, which we'll denote by
      $H^{2n+1}$ or $H$, is a connected, simply connected real nilpotent
      Lie group. It is modelled on $\mathbb{R}^{2n+1}$ with the group law
      \[ (x,y,w) \cdot (x',y',w') = (x+x', y+y', w+w'+\tfrac{1}{2}
      (x \cdot y' - y \cdot x')) \]
      for $x, y \in \mathbb{R}^n, w \in \mathbb{R}$.

      Its Lie algebra $\mathfrak{h}$ has basis $\{ X_1, \dots, X_n,
      Y_1, \dots, Y_n, W \}$ and non-zero commutation relations
      $[X_j, Y_j] = W = - [Y_j, X_j]$.

      Let $X_j$ also denote the left-invariant vector field on $H$
      given by left translation of $X_j \in \mathfrak{h}$,
      which we identify with an element of the tangent space at the
      identity.

      We define complex vector fields $Z_j, Z_{\bar{j}}$ on $H$ by
      $Z_j := \tfrac{1}{\sqrt{2}} (X_j - i Y_j),
      Z_{\bar{j}} := \tfrac{1}{\sqrt{2}} (X_j + i Y_j).$
      Alternatively, with the same definitions, we consider
      $Z_j$ and $Z_{\bar{j}}$ to be elements of $u(\mathfrak{h})$,
      the universal enveloping algebra of $\mathfrak{h}$.

      We choose a left-invariant metric on $H$ such that
      $\{ X_j, Y_j, W \}$ is an orthonormal basis for $T_p H$ at each point
      $p$ of $H$. Then $\{ Z_1, \dots, Z_n, Z_{\bar{1}}, \dots,
      Z_{\bar{n}}, W \}$ is an orthonormal basis for the complexified
      tangent space at each point.

      Let $\{ \tau^1, \dots, \tau^n, \tau^{\bar{1}}, \dots, \tau^{\bar{n}},
      \tau^w \}$ be the basis of 1-forms dual to
      $\{ Z_1, \dots, Z_n,$ $Z_{\bar{1}}, \dots, Z_{\bar{n}}, W \}$.

We define $\hat{G}$ to be the set of (unitary) equivalence classes 
of irreducible unitary representations of a locally compact group $G$.

If $\pi$ is a unitary representation of a group $G$ on a 
Hilbert space
${\mathcal H}_\pi$, then there is an induced representation of $L^1(G)$
on ${\mathcal H}_\pi$. That is, take any element $f$ in $L^1(G)$:
we define
\[ \pi(f) := \int_G f(x) \pi(x) dx. \]
If we have any elements $u, v \in {\mathcal H}_\pi$, then
\[ \innprod{\pi(f)u}{v} = \int_G f(x) \innprod{\pi(x)u}{v} dx. \]
The operator $\pi(f)$ is known as the \emph{group Fourier transform}
of $f$ (see \cite{dixmierc} or \cite{folland}).

Let $\mathcal{J}^{1}$ be $L^{1}(G) \cap L^{2}(G)$ and
$\mathcal{J}^{2}$ be the set of finite linear combinations of
elements of the form $f * g$, for $f,g \in \mathcal{J}^{1}$.
(Note the similarities to Hilbert-Schmidt and nuclear or trace-class
operators.)

The following theorem, the abstract Plancherel theorem,
was first proved in
\cite{mautner}, \cite{segal}; the formulation below is taken from
\cite{dixmierc} and \cite{follandnew}.

\begin{theorem} \label{T:AbPlThm}
Let $G$ be a type I, unimodular, separable, locally compact group.
Then there exists a measurable field of irreducible representations
$\pi_\zeta$ over $\hat{G}$ such that $\pi_\zeta$ belongs to
the equivalence class $\zeta$. We identify $\pi_\zeta$ with $\zeta$,
and write ${\mathcal H}_\zeta$ for the Hilbert space which
$\pi_\zeta$ acts on. Let $t_\zeta$ be the trace $T \otimes 1
\mapsto {\rm Tr}_\zeta (T)$ (the Hilbert-Schmidt trace)
on the positive operators in
$B({\mathcal H}_\zeta) \otimes \mathbb{C}$.

Let $\pi_L$ and $\pi_R$ be the left and right regular representations
of $G$, and let ${\mathcal U}$ and $\mathcal V$ be the von Neumann
algebras on $L^2(G)$ generated by $\pi_L(G)$ and $\pi_R(G)$.
      Let $t$ be the trace on ${\mathcal U}^+$
	defined as above.

Then there exists a positive measure $\mu$ on $\hat{G}$ and an
isomorphism $W$ from $L^2(G)$ to $\int_{\hat{G}}^\oplus
({\mathcal H}_\zeta \otimes \overline{\mathcal H}_\zeta)
\, d \mu (\zeta)$ such that:
\begin{enumerate}
\item
$W$ transforms
$\pi_L$ into $\int_{\hat{G}}^\oplus (\zeta \otimes 1) \, d\mu(\zeta)$,
$\pi_R$ into $\int_{\hat{G}}^\oplus (1 \otimes \bar{\zeta}) \, d\mu(\zeta)$,
$\mathcal U$ into $\int_{\hat{G}}^\oplus (B({\mathcal H}_\zeta)
\otimes \mathbb{C}) \, d\mu(\zeta)$,
$\mathcal V$ into $\int_{\hat{G}}^\oplus (\mathbb{C} \otimes B
(\overline{\mathcal H}_\zeta)) \, d\mu(\zeta)$, and
$t$ into $\int_{\hat{G}}^\oplus t_\zeta \, d\mu(\zeta)$.
\item
If $h \in \mathcal{J}^{2}$ and $x \in G$, then we have the Fourier inversion
formula for $G$:
      \begin{equation} \label{E:FIT4G}
      h(x) = \int_{\hat{G}} {\rm Tr} (\zeta(x) \zeta(h)) \, d \mu (\zeta).
      \end{equation}
In particular, if $u \in L^1(G) \cap L^2(G) = \mathcal{J}^{1}$, we have
\[ \int_G \abs{u(s)}^2 \, ds = \int_{\hat{G}} {\rm Tr} (\zeta(u) \zeta(u)^*)
\, d \mu (\zeta), \]
the Plancherel formula for $G$.
\end{enumerate}
\end{theorem}

Note that we write $\int_{\hat{G}}^\oplus (\zeta \otimes 1) \, d\mu(\zeta)$
rather than $\int_{\hat{G}}^\oplus \zeta \, d\mu(\zeta)$ and 
$\int_{\hat{G}}^\oplus (B({\mathcal H}_\zeta)
\otimes \mathbb{C}) \, d\mu(\zeta)$ rather than 
$\int_{\hat{G}}^\oplus B({\mathcal H}_\zeta) \, d\mu(\zeta)$;
this is to clarify the action of these operators on 
$\int_{\hat{G}}^\oplus ({\mathcal H}_\zeta \otimes \overline{\mathcal H}_\zeta)
\, d \mu (\zeta)$.

The measure $\mu$ is known as the Plancherel measure of $\hat{G}$
(associated with the Haar measure of $G$).

For the Heisenberg group, the Plancherel measure $\mu$ is zero except on
representations $\beta_k$,
where $\beta_k$ is the Fock-Bargmann representation of $H$ with
parameter $k$ (for $k \in {\mathbb R}^*$). This representation is
irreducible and acts on the Fock space $\mathcal{F}^k_n$
defined by
      \[ \mathcal{F}^k_n = \{ F : F \, \text{is entire on $\mathbb{C}^n$ and}
      \, \int_{\mathbb{C}^n} \abs{F(z)}^2 e^{-k z \cdot \bar{z}/4} dz <
      \infty \}. \]

In fact, we're more interested in the conjugate representation
$\bar{\beta}_k$, which acts on the anti-Fock space $\mathcal{F}^{-k}_n$,
where $F \in \mathcal{F}^{-k}_n$ iff $\bar{F} \in \mathcal{F}^{k}_n$.
This representation is defined by
      \[ \bar{\beta}_{k}(p,q,w)F(\bar{z})=
 e^{-ikw - \tfrac{1}{4}k(p^2+q^2)-\tfrac{1}{2}k\bar{z} \cdot (p+iq)}
 F(\bar{z}+p-iq). \]

%(1-dimensional) center $Z$ of $H$ it
%satisfies $$\beta_k( {\rm exp} (tW)) = e^{ikt} \, \text{Id},$$
%where $t \in R$ (so that $tW \in \mathfrak{z}$, the centre of $\mathfrak{h}$,
%and ${\rm exp}(tW) \in Z$).

      Thus, the Plancherel theorem for $H$ implies that
      \begin{equation} \label{E:PL4H}
      L^2 (H) \isom \int^\oplus_{k \in \mathbb{R}}
      \mathcal{F}^k_n \otimes \mathcal{F}^{-k}_n \abs{k}^n dk.
      \end{equation}
      Under this decomposition, the right regular representation
      $\pi_R$ of $H$ on $L^2(H)$ is given by
      \[ \pi_R = \int^\oplus_{k \in \mathbb{R}} (\id \otimes
      \bar{\beta}_k) \abs{k}^n dk.\]
 
      From the representation 
$\bar{\beta}_k$ of $H$, we have a representation (also denoted by
      $\bar{\beta}_k$) of $u(\mathfrak{h})$ on the $C^\infty$ vectors
      of $\mathcal{F}^{-k}_n$ (see \cite{c&g, garding}), given by
      \[ \bar{\beta}_k(Z_j) = - \tfrac{1}{\sqrt{2}} k \bar{z}_j, \,
      \bar{\beta}_k (Z_{\bar{j}}) = \sqrt{2} \partial_{\bar{z}_j}, \,
      \bar{\beta}_k (W) = - ik.
      \]

      For any multi-index $\beta \in \mathbb{Z}_+^n$, we define a function
      $\psi_\beta(k)$ by
      \[ \psi_\beta (k) := \left(\frac{k}{2\pi}\right)^{n/2} \left({\frac{ik}{2}}
 \right)^{\abs{\beta}/2} \frac{\bar{z}^\beta}{\sqrt{\beta!}}. \]
      Then the set $\{ \psi_\beta(k) : \beta \in \mathbb{Z}_+^n \}$ is
      a complete orthonormal basis of $\mathcal{F}^{-k}_n$ 
(see \cite{folland}).

      The action of the above operators on this basis is given by
      \begin{eqnarray*}
      \bar{\beta}_k (Z_j) (\psi_\beta(k)) & = & -i \sqrt{k}
        \sqrt{\beta_j + 1} \, \psi_{\beta + e_j} (k), \\
      \bar{\beta}_k (Z_{\bar{j}}) (\psi_\beta(k)) & = & - i \sqrt{k}
        \sqrt{\beta_j} \, \psi_{\beta- e_j}(k)
      \end{eqnarray*}
      where $e_j$ is the multi-index with $1$ in the $j$th place and
      zeros elsewhere.

We define creation and annihilation operators $a_j, a_j^*$ which act on
$\mathcal{F}^{-k}_n$. Let $a_j$ be the operator $i k^{-1/2} \bar{\beta_k}
(Z_{\bar{j}})$, and $a_j^*$ the operator $i k^{-1/2} \bar{\beta_k} (Z_j)$.
Then $[ a_j, a_j^*] = \id$.
We call $a_j^*$ a creation operator and $a_j$ an annihilation operator.
Note that 
\begin{eqnarray*}
      a_j^* \psi_\beta(k) & = & 
	\sqrt{\beta_j + 1} \, \psi_{\beta + e_j} (k), \\
      a_j \psi_\beta(k) & = & 
        \sqrt{\beta_j} \, \psi_{\beta- e_j}(k)
      \end{eqnarray*}

\section{An explicit formula for the Laplacian}

In this section, we begin to explicitly analyse the action of the Laplacian.

For $d: \Lambda^p_{(2)}H \otimes {\mathbb C} \rightarrow \Lambda^{p+1}_{(2)}H
\otimes {\mathbb C}$ the (complexified) exterior derivative on $L^2$ $p$-forms
and
$d^*: \Lambda^p_{(2)}H \otimes {\mathbb C} \rightarrow \Lambda^{p-1}_{(2)}H
\otimes {\mathbb C}$ its adjoint, the Laplacian on $p$-forms is
defined to be
\[ \triangle = dd^* + d^* d : \Lambda^p_{(2)}H \otimes {\mathbb C} \rightarrow
\Lambda^p_{(2)}H \otimes {\mathbb C}. \]
It will also be denoted by $\triangle_p$ or $\triangle_{p,n}$, when the
degree of the forms and/or the dimension of the group that the Laplacian
is acting on is important.
      Note that the domain of the Laplacian is the first Sobolev space of
      $p$-forms; since this is dense in the space of $L^2$ $p$-forms, we
      assume
      for the purposes of this
discussion that the Laplacian acts on $L^2$ $p$-forms.
(For more on this, see \cite{atiyah, dodziuk}.)

      In particular, the Laplacian on functions is given by
      \[
      \triangle_{0,n} = \sum_{j=1}^n (-Z_j Z_{\bar{j}} -Z_{\bar{j}}
      Z_j ) -W^2, \]
which implies that, acting on
$\mathcal{F}^{-k}_{n}$, $\triangle_{0,n} (k) = \sum_{j=1}^n (2k a_j^* a_j
      +k) + k^2$.
      In particular, on the basis elements, 
$\triangle_{0,n} (k) \psi_{\beta} (k) = (2k \abs{\beta}
      +nk + k^2) \psi_{\beta} (k)$ where $\abs{\beta} = \beta_1
      + \ldots + \beta_n$.

      By inspection, the lowest eigenvalue of $\triangle_{0,n} (k)$ is $nk +
      k^2$.
%as is well-known.
%see pesce etc.
      Furthermore, the eigenvalue corresponding to $\psi_{\beta}
      (k)$ depends only on $\abs{\beta}$, and not on any other function of
      $\beta$.

We begin by calculating explicitly the form of $d$ and $d^*$ acting on
$p$-forms.

\begin{lemma}
The actions of $d$ and $d^*$ on $p$-forms on
$H^{2n+1}$ are given by
\begin{eqnarray*}
d & = & \left( \sum_{j=1}^n e(\tau^j)Z_j + e(\tau^{\bar{j}})Z_{\bar{j}}
        \right)
		+ e(\tau^w)W - i \sum_{j=1}^n e(\tau^j)e(\tau^{\bar{j}})i(W) \\
d^* & = & -\left( \sum_{j=1}^n i(Z_{\bar{j}})Z_j + i(Z_j)Z_{\bar{j}}
          \right) - i(W)W
	      + i\sum_{j=1}^n e(\tau^w)i(Z_{\bar{j}})i(Z_j)
\end{eqnarray*}
where $e(\tau)$ denotes exterior multiplication by the $1$-form $\tau$
and $i(V)$ denotes contraction by the vector field $V$.
\end{lemma}

The proof of this lemma uses the Leibnitz rule (giving the first few
terms
in the above formula for $d$, which are the same as those for $d$ on
functions) and the fact that for any $1$-form $\eta$ and vector fields
$X,Y$,
$d \eta (X,Y) = X \eta (Y) - Y \eta (X) - 1/2 \eta ([X,Y])$ (see for 
example \cite{spivak}).
The action for $d$ on $2$-forms is unremarkable since the Heisenberg
group is a $2$-step nilpotent Lie group.

Using these formulae for $d$ and $d^*$, we can explicitly calculate
the form of $\triangle_{p,n}$ (again in terms of $e(*)$ and $i(*)$'s).
(The details of this calculation are given in an appendix.)
Here we write $\triangle_{p,n}$ (and $\triangle_{p,n}(k)$) as a matrix,
considering a $p$-form to be a
$\tbinom{2n+1}{p}$ vector
- again using the triviality of the tangent bundle of $H^{2n+1}$.
Recall that for any operator $A$ acting on $\Lambda^*(H) \otimes {\mathbb
C}$ or on $L^2(H) \otimes {\mathbb C}$, we denote the decomposition in
the representation $\bar{\beta}(k)$ by $A(k)$.

The Laplacian on $p$-forms, acting on $H^{2n+1}$, is given by:
\begin{align} \label{E:full_lap}
\triangle_{p,n} & = -W^2 + \sum_{j=1}^n \left( -2Z_jZ_{\bar{j}} 
		+ iW (i(Z_j)e(\tau^j) +
             e(\tau^{\bj})i(Z_{\bj})) \right.\nn \\
          & + ie(\tau^w)(i(Z_{\bar{j}})Z_j
           - i(Z_j)Z_{\bj}) % \nn \\
          - i (e(\tau^j)Z_j - e(\tau^{\bj})Z_{\bj})i(W) \nn \\
          & + \sum_{k=1,k \neq j}^n
          e(\tau^j)e(\tau^{\bj})i(Z_{\bk})i(Z_k) \nn \\
          & \left. + e(\tau^j)i(Z_j)e(\tau^{\bj})i(Z_{\bj})i(W)e(\tau^w) %\\
          + i(Z_j)e(\tau^j)i(Z_{\bj})e(\tau^{\bj})e(\tau^w)i(W) \right)
\end{align} 
(This formula is derived in Appendix A.)

After the transform corresponding to the conjugate
Fock-Bargmann representation
with parameter $k$, this operator becomes:
\begin{align} \label{E:lap_rep}
\triangle_{p,n}(k)
          & = k^2 + \sum_{j=1}^n \left( 2ka_j^*a_j + k i(Z_j)e(\tau^j) +
            k e(\tau^{\bj})i(Z_{\bj}) \right. \nn \\
          & + \sqrt{k}e(\tau^w)(i(Z_{\bar{j}})a_j^*
           - i(Z_j)a_j) %\\
           + \sqrt{k} i(W)(e(\tau^j)a_j^*
           - e(\tau^{\bj})a_j)  \nn \\
          & + \sum_{k=1,k \neq j}^n e(\tau^j)e(\tau^{\bj})i(Z_{\bk})
           i(Z_k) \nn \\
          & \left. + e(\tau^j)i(Z_j)e(\tau^{\bj})i(Z_{\bj})i(W)e(\tau^w) %\\
           + i(Z_j)e(\tau^j)i(Z_{\bj})e(\tau^{\bj})e(\tau^w)i(W) \right)
\end{align}
Using this last formula, we could explicitly calculate all the eigenvalues of
$\triangle_{p,n}(k)$
for certain (small) values of $n$ and $p$, writing the Laplacian
globally
as a matrix (since the tangent space of $H^{2n+1}$ is trivial). However, the
size of this matrix is $\tbinom{2n+1}{p}$, as implied above,
and so will grow polynomially as $n$ and $p$ increase.

      We note instead that we can define the
      following operators.
      \begin{definition} \label{D:theta}
      For $j=1, \ldots, n$, we define $\theta_{j}$ to be a map from
      $\Lambda^{p}_{(2)}(H) \otimes \mathbb{C}$ and $\theta_{j}^{*}$
      to be its adjoint, given by the following formulae:
      \begin{eqnarray*}
      \theta_j & = & e(\tau^j)Z_j + e(\tau^{\bar{j}})Z_{\bar{j}}
                -ie(\tau^j)e(\tau^{\bar{j}})i(W) \\
      \theta_j^* & = & -i(Z_j)Z_{\bar{j}} - i(Z_{\bar{j}})Z_j
                +ie(\tau^w)i(Z_{\bar{j}})i(Z_j).
      \end{eqnarray*}
      \end{definition}
      We can then rewrite $d$ and $d^{*}$ as 
       $d= \sum_j \theta_j + e(\tau^w)W$ and $d^* = \sum_j \theta_j^*
      - i(W)W$.

	Writing $\triangle_{p,n}(k)$ in terms of the operators $\theta_j(k),
	\theta_j^*(k), e(\tau^w)k$ and $i(W)k$ gives us further information
	about the spectrum of $\triangle_{p,n}(k)$; in particular, we find
	a lower bound on the spectrum for all $p$ and $n$, which is achieved
	for $p=n$.
	
      \begin{lemma} \label{L:lower_bound}
      The operator $\triangle_{p,n}(k)$ satisfies the inequality:
      \[ \triangle_{p,n}(k) \ge k^2 \id. \]
      In particular, $\triangle_{n,n}(k)$ has lowest eigenvalue
      $k^2$.
      \end{lemma}

      \begin{proof}
Since $e(\tau^w)$ and $\theta_j^*$ anticommute, as do $i(W)$ and 
$\theta_j$, for all $j$, we have that 
      \begin{eqnarray*}
      \triangle & = & (\sum_j \theta_j)(\sum_m \theta_m)^* + (\sum_m \theta_m)^*
      (\sum_j \theta_j) - W^2 \\
\implies      \triangle(k) & \ge & k^2
      \end{eqnarray*}

      This is a lower bound on the eigenvalues of $\triangle_{p,n} (k)$
      for all $n$ and $p$. However, if $n=p$, we know (from \cite{lott}) that
      there is an eigenvector $v$ of $\triangle_{n,n}(k)$, 
      \[ v:= f \tau^1 \wedge \ldots \wedge \tau^n, \]
      where $f \in \ker Z_{\bar{1}}(k) \cap \ldots \cap \ker Z_{\bar{n}}(k)$,
      such that
      $ \triangle_{n,n} (k) v = k^2 v;$
      thus $k^2$ is in fact the lowest eigenvalue of $\triangle_{n,n}
      (k)$ for all $n$.
      \end{proof}

      \section{Commuting operators}

     In this section, we define a partition of %the Hilbert space
$\mathcal{F}^{-k}_{n} \otimes \Lambda^{p}(\mathfrak{h}^{*})$
into subspaces, using a collection of commuting operators.

For $j=1,\dots,n$, we define $U_{jj}$ to be the operator on
$\mathcal{F}^{-k}_{n} \otimes \Lambda^{p}(\mathfrak{h}^{*})$
given by
\[ U_{jj}:= a_{j}^{*} a_{j} - e(\tau^{j})i(Z_{j}) + e(\tau^{\bj})
i(Z_{\bj}).\]
It should be clear from this definition that $[U_{jj},U_{ll}]=0$
for all $j \neq l$, and that this operator is self-adjoint:
$U_{jj}^*=U_{jj}$.

Define the set $S:= \{ \gamma \in \mathbb{Z}^{n} : \gamma_{j} \ge
      -1, j=1,\dots,n,$ and at most $p$ of the indices $\gamma_{j}$
     are equal to $-1 \}$.
For multi-indices $\gamma$ in $S$,
we define the subspace $V^{p,n,\gamma}$ to be the simultaneous
eigenspace of $U_{11}, \dots, U_{nn}$, with eigenvalues
$\gamma_{1}, \dots, \gamma_{n}$.       (We sometimes omit the mention of $n$.)
That is, if we write $E_{\lambda} A$ for
the eigenspace of an operator $A$ corresponding to the eigenvalue
$\lambda$, then $V^{p,n,\gamma}$ is given by
\[ V^{p,n,\gamma} := E_{\gamma_1} U_{11} \cap \dots E_{\gamma_{n}}
U_{nn} \cap (\mathcal{F}^{-k}_{n} \otimes \Lambda^{p}(\mathfrak{h}^{*})). \]

For example, the $p$-form $\psi_{\gamma+I-J} \tau^{I} \wedge
\tau^{\bar{J}}$ is in $V^{p,n,\gamma}$, where $I$ and $J$ are both multi-indices,
with entries either 0 or 1, $\abs{I}+\abs{J}=p$
      and if $I=e_{i_1} + \ldots + e_{i_m}$, then $\tau^I = \tau^{i_1}
      \wedge \ldots \wedge \tau^{i_m}$ (and similarly for
      $\tau^{\bar{J}}$).

      Note that we can have $\gamma_j=-1$ for some $j=1,\dots,n$,
      but this means that
      every element of $V^{p,n,\gamma}$ would have to be of the form
      $\tau^j \wedge v$ for some $v \in V^{p-1,n,\gamma+e_j}$;
      thus at most $p$ of the $\gamma_j$'s can be -1. The remainder of the
      indices of $\gamma$ must be non-negative. 

      In fact, the collection of the subspaces
      $V^{p, n, \gamma}$ for all values of
      $\gamma$ in $S$ is a partition: 
      \[ \mathcal{F}^{-k}_{n} \otimes \Lambda^{p} (\mathfrak{h}^{*})
      = \oplus_{\gamma \in S} V^{p,n,\gamma}. \]

      The subspace $V^{0,n,\gamma}$ consists of (complex) scalar
      multiples of $\psi_{\gamma} (k)$;
      the subspace $V^{p,n, \gamma}$ also corresponds to $\psi_{\gamma}
      (k)$ in some sense, but with dimension $\tbinom{2n+1}{p}$.
  
      The usefulness of this definition is due to the following theorem.

\begin{theorem}
Let $d(k)$ and $d^*(k)$ represent the exterior differential
and its adjoint respectively in the representation $\bar{\beta}_k$. Then
$d(k)$ maps $V^{p,n,\gamma}$ to $V^{p+1,n,\gamma}$ for $p < 2n+1$, and
$d^*(k)$ maps $V^{p,n,\gamma}$ to $V^{p-1,n,\gamma}$, for $p \ge 1$.
So $V^{p,n, \gamma}$ is a $\triangle_{p,n}(k)$-invariant subspace of
      $\mathcal{F}^{-k}_n \otimes \Lambda^p (\mathfrak{h}^*)$.
\end{theorem}

\begin{proof}
We prove that $[U_{jj},\theta_{j}(k)]=0$, and thus that
$[U_{jj},d(k)]=0 = [U_{jj},d^{*}(k)]$ for all $j$, which means that
$[U_{jj},\triangle_{p,n}(k)]=0$ for all $j$.

\begin{eqnarray*}
[\theta_j(k), U_{jj}]
& = & \sqrt{-1} ([k^{-1/2}e(\tau^j)a_j^*, -e(\tau^j)i(Z_j)]
	+ [k^{-1/2}e(\tau^{\bar{j}}) a_j, a_j^* a_j] \\
& & + [-e(\tau^j)e(\tau^{\bar{j}})i(W),-e(\tau^j)i(Z_j)] \\
& & + [k^{-1/2}e(\tau^j)a_j^*, a_j^* a_j]
	+ [k^{-1/2}e(\tau^{\bar{j}})a_j, e(\tau^{\bar{j}}) i(Z_{\bar{j}})]
	\\
& &+ [-e(\tau^j)e(\tau^{\bar{j}})i(W),e(\tau^{\bar{j}}) i(Z_{\bar{j}}]) \\
& = & \sqrt{-1}(k^{-1/2}e(\tau^j)a_j^* + k^{-1/2}e(\tau^{\bar{i}})a_j
	-e(\tau^j)e(\tau^{\bar{j}})i(W) \\
& & - k^{-1/2} e(\tau^j)a_j^* - k^{-1/2} e(\tau^{\bar{j}})a_j
	-e(\tau^{\bar{j}})e(\tau^j)i(W) )\\
& = & 0.
\end{eqnarray*}
Clearly, $U_{jj}$ also commutes with $\theta_{l}(k)$ (for $l \neq j$),
since different operators are involved.
So $U_{jj}$ commutes with $d(k)$; then since $U_{jj}^{*}=U_{jj}$,
this means that $U_{jj}$ also commutes with $d^{*}(k)$, and thus that
$U_{jj}$ commutes with $\triangle_{p,n}(k)$.

This means that any eigenspace of $U_{jj}$ will be preserved by
$\triangle_{p,n}(k)$.
But this is true for all $j$, so the subspace $V^{p,n,\gamma}$ is
$\triangle_{p,n}(k)$-invariant.
\end{proof}

      We can now study the eigenvalues of $\triangle_{p,n} (k)$
      restricted to $V^{p, \gamma}$. In fact, we'll also be interested in even
smaller subspaces. For this, the following definition will be useful.

\begin{definition}
Fix $k>0$. Let $V$ be a subspace of
$\mathcal{F}^{-k}_{n} \otimes \Lambda^p
((\mathfrak{h}^{2n+1})^*)$
and let $W$ be a subspace of
$\mathcal{F}^{-k}_n \otimes \Lambda^q ((\mathfrak{h}^{2m+1})^*)$
for some $n,m,p,q$ (such that $p \le 2n+1$ and $q \le 2m+1$), so that
$\triangle_{p,n}(k)$ acts on $V$ and $\triangle_{q,m}(k)$ acts on $W$.
Suppose also that $V$ is $\triangle_{p,n}(k)$-invariant and $W$ is
$\triangle_{q,m}(k)$-invariant.
Then we say that $V$ and $W$ are {\bf spectrally equivalent} if
$\triangle_{p,n}(k)$ acting on $V$ has the same eigenvalues (including multiplicity)
as $\triangle_{q,m}(k)$ acting on $W$.
\end{definition}

\begin{remark}
This is true if and only if there is a linear isomorphism
$j$ from $V$ to $W$
which commutes with $\triangle(k)$, i.e. such that
$j \triangle_{p,n}(k) = \triangle_{q,m}(k) j.$
From either of these conditions, we can see that spectral equivalence is
indeed an equivalence relation.
\end{remark}

We now introduce an operator on $p$-forms, as a first step in calculating the
eigenvalues of $\triangle_{p,n}(k)$.

\begin{definition}
The $(1,2)$ transposition operator is an operator on
$\mathcal{F}^{-k}_n \otimes \Lambda^p(\mathfrak{h}^*)$,
denoted by $U_{12}$ and defined to be
      $$ U_{12} := a_1^* a_2 -e(\tau^2) i(Z_1) + e(\tau^{\bar{1}})
      i(Z_{\bar{2}}).$$
Similarly, we define $U_{ij}$, the $(i,j)$ transposition operator,
(for $i \neq j, i,j =1, \ldots, n$) to be the
operator given by
$$ U_{ij} := a_i^* a_j -e(\tau^j) i(Z_i) + e(\tau^{\bar{i}})
      i(Z_{\bar{j}}).$$
\end{definition}

      From the above definition, we see that
      the $(j,i)$ transposition operator $U_{ji}$ is the
      adjoint of $U_{ij}$.
Also, $U_{ij}$ is ``usually'' an isomorphism, as proved in the following lemma.

\begin{lemma} \label{L:Uker}
\begin{enumerate}
\item If $v \in \ker U_{ij} \cap V^{p,n,\gamma}$ for
some $\gamma$, then $\gamma_j = -1,0$ or $1$. %\\
\item If $\gamma=(\gamma_1,\gamma_2,\ldots, \gamma_n)$, $\gamma_i \ge 1$ and
$\gamma_j \ge 2$ (for $i \neq j$), then the restriction of $U_{ij}$ to
$V^{p,n,\gamma}$
is a linear isomorphism between $V^{p,n,\gamma}$ and
$V^{p,n,\gamma+e_i-e_j}$.
\end{enumerate}
\end{lemma}

\begin{proof}
The proof of (i) can be found in the appendix.
 To prove (ii), we note from (i) that $U_{ij}$ restricted to
 $V^{p,n,\gamma}$
 is 1-1 (since $\gamma_j \ge 2$). Now the orthogonal complement of the image of
 $U_{ij}$ is
 the kernel of the adjoint map.
 But the adjoint of $U_{ij}$ on $V^{p,n,\gamma}$ is $U_{ji}$ restricted
 to $V^{p,n,\gamma+e_i-e_j}$, which has kernel $\{0\}$, again by (i)
 (since $\gamma_i+1 \ge 2$). So
 this map is onto and thus an isomorphism.
 \end{proof}

 \begin{remark}
 Note that $U_{ij} U_{ji}$ is not the identity; however, it
 is an automorphism on $V^{p,n,\gamma}$ for ``generic'' $\gamma$,
 and since it commutes with $\triangle(k)$, it preserves eigenspaces.
 \end{remark}

 We also have:

 \begin{lemma}\label{L:Ucomm}
 \begin{enumerate}
 \item \label{L:Ud} For all $i,j$ and $p$,
 $[d(k),U_{ij}]=0$; also $[d^*(k),U_{ij}]=0$ and thus $[\triangle_p(k),U_{ij}]=0$.
 That is, the $(i,j)$ transposition operator commutes with the Laplacian on $p$-forms
 in the representation $\bar{\beta}_k$. %\\
 \item If $\gamma=(\gamma_1,\gamma_2,\ldots, \gamma_n)$, $\gamma_i \ge 1$ and
$\gamma_j \ge 2$, then $V^{p,n,\gamma}$ and $V^{p,n,\gamma+e_i-e_j}$
are spectrally equivalent. %\\
 \item \label{L:gamma_beta}
 For any $\beta$, $\gamma$ multi-indices such that $\abs{\beta}=\abs{\gamma}$
 and $\beta_i \ge 1$, $\gamma_i \ge 1$ for all $i=1,2,\ldots,n$, the
 subspaces
 $V^{p,n,\beta}$ and $V^{p,n,\gamma}$ are spectrally equivalent.
 \end{enumerate}
 \end{lemma}

\begin{proof}
In proving (i), note that for any $l \neq i, l \neq j$, we have
that $[\theta_l(k), U_{ij}]=0$, so that we only need prove that
$[\theta_i(k)+\theta_j(k), U_{ij}]=0$.
Now
\begin{eqnarray*}
& & [\theta_i(k) + \theta_j(k), U_{ij}] \\
& = & \sqrt{-1} ([k^{-1/2}e(\tau^i)a_i^*, -e(\tau^j)i(Z_i)]
	+ [k^{-1/2}e(\tau^{\bar{i}}) a_i, a_i^* a_j] \\
& &	+ [-e(\tau^i)e(\tau^{\bar{i}})i(W),-e(\tau^j)i(Z_i)] \\
& & + [k^{-1/2}e(\tau^j)a_j^*, a_i^* a_j]
	+ [k^{-1/2}e(\tau^{\bar{j}})a_j, e(\tau^{\bar{i}}) i(Z_{\bar{j}})]
	\\
& &	+ [-e(\tau^j)e(\tau^{\bar{j}})i(W),e(\tau^{\bar{i}}) i(Z_{\bar{j}}]) \\
& = & \sqrt{-1}(k^{-1/2}e(\tau^j)a_i^* + k^{-1/2}e(\tau^{\bar{i}})a_j
	-e(\tau^j)e(\tau^{\bar{i}})i(W) \\
& & - k^{-1/2} e(\tau^j)a_i^* - k^{-1/2} e(\tau^{\bar{i}})a_j
	-e(\tau^{\bar{i}})e(\tau^j)i(W) )\\
& = & 0.
\end{eqnarray*}
One can prove that $[d^*(k),U_{ij}]=0$ in a similar way, or take the
adjoint of the equation $[d(k),U_{ji}]=0$. It then follows that
$[\triangle(k),U_{ij}]=0$.

     (ii) follows from (i) and from Lemma \ref{L:Uker} (i); (iii)
     is easily proved by repeated use of (ii) for selected values of
     $i$ and $j$.
\end{proof}

      Thus the eigenvalues of $\triangle_{p,n} (k)$ on different
      $V^{p, \gamma}$ also depend only on $| \gamma |$ for generic
      $\gamma$, as is the case when $p=0$ for any $\gamma$.

%      \begin{remark}
%      What is the relevance of this operator? We provide some
%      motivation for its definition, beginning with the situation of $L^{2}$
%      functions on $H$.

 %     The operator $a_1^* a_2$
 %     acting on $\psi_{\beta}(k)$ gives a multiple of $\psi_{\beta+
 %     e_1 - e_2}(k)$; i.e. acting on an eigenvector of $\triangle_0
 %     (k)$, this operator gives another eigenvector with the same
 %     eigenvalue. This statement is equivalent to the equation
 %     $$[ \triangle_0 (k), a_1^* a_2] = 0. $$

%      It may be helpful at this point to present the quantum mechanical
%      interpretation of this equation.
%      The creation operator $a_j^*$ represents the creation of a
%      (bosonic) particle
%      in the $j$th state, and the annihilation operator $a_j$ the
%      destruction of such a particle; the Laplacian on functions,
%      $\triangle_0(k)$, is the number operator (counting the total
%      number of particles in all states) plus a potential. The operator
%      $a_1^* a_2$, then, represents changing a particle from state $2$
%      to state $1$, which leaves the total number of particles unchanged.

%Generalising this operator to $p$-forms gives $U_{12}$, where the
%extra terms are added so that $U_{12}$ commutes with the Laplacian
%on $p$-forms. The quantum
%mechanical interpretation of this transposition operator is less clear,
%but may also involve interchanges
%of fermionic particles.
%\end{remark}

The following definitions rely on the fact that
$\mathfrak{h}_{2n+1}$ is symmetric with respect
to the basis elements $X_1, Y_1, X_2, Y_2, \ldots, X_n, Y_n$: i.e.,
if $X_2$ and $Y_2$ are interchanged with $X_1$ and $Y_1$, then
the commutation relations are unchanged.

\begin{definition}
 We define an action of $S_n$ (the permutation group on
$n$ symbols) on $\mathbb{Z}^n$ by:
$$ \sigma \cdot (\beta_1,\beta_2,\ldots,\beta_n) = (\beta_{\sigma(1)},
\beta_{\sigma(2)},\ldots,\beta_{\sigma(n)})$$
for $\sigma \in S_n$ and $\beta_1,\ldots,\beta_n \in \mathbb{Z}$.
For example, for $\beta=(\beta_1,\ldots,\beta_n) \in \mathbb{Z}^n$, we
have $(12) \cdot \beta = (\beta_2,\beta_1,\ldots,\beta_n)$.
\end{definition}

We next define an isometry for all pairs of distinct numbers $i,j$,
utilising the symmetry of the Lie algebra of the Heisenberg group.

\begin{definition} Let $\chi_{ij}$ be the Lie algebra isomorphism on
$\mathfrak{h}^{2n+1}$ (for $i \neq j, i,j \le n$) defined by linearity
and action
$\chi_{ij}: Z_i \mapsto Z_j, Z_j \mapsto Z_i, Z_{\bar{i}} \mapsto
Z_{\bar{j}}, Z_{\bar{j}} \mapsto Z_{\bar{i}},$ and $V \mapsto V$ if $V$
is orthogonal to $Z_i, Z_j, Z_{\bar{i}}$ and $Z_{\bar{j}}$.
\end{definition}
This map $\chi_{ij}$ is an isometry with respect to the inner product
that we have chosen on $\mathfrak{h}^{2n+1}$.

It induces a map $\tilde{\chi}_{ij}$ on $L^{2}$ $p$-forms on $H$, and
this factors through the representation to give a map on $\mathcal{F}^{-k}_{n}
\otimes
\Lambda^{p} (\mathfrak{h}^{*})$, which we'll
also denote by $\chi_{ij}$.
This map is linear, multiplicative, and has action
\begin{eqnarray*}
\tau^{i} \mapsto \tau^{j}; & \tau^{j} \mapsto \tau^{i}; & \tau^{m}
\mapsto \tau^{m} \, \text{if} \, m \neq i,j; \\
\tau^{\bar{i}} \mapsto \tau^{\bar{j}}; & \tau^{\bar{j}} \mapsto
\tau^{\bar{i}}; & \tau^{\bar{m}} \mapsto \tau^{\bar{m}} \, \text{if}
\, m \neq i,j; \\
\tau^{w} \mapsto \tau^{w}; & \psi_{\beta} (k) \mapsto
\psi_{(ij) \cdot \beta} (k) &
\end{eqnarray*}

So the operator $\chi_{ij}$ is an isometry from
$V^{p,n,\gamma}$ to $V^{p,n,(ij) \cdot \gamma}$
and commutes with the Laplacian $\triangle_{p,n}(k)$.

We call $\chi_{ij}$ the $(i,j)$ symmetry operator, or simply a symmetry
operator.

These operators can be
used to prove that the eigenvalues of the Laplacian in the
representation $\bar{\beta}_k$ on the subspace $V^{p,n,\gamma}$
are symmetric in the entries of $\gamma$. That is, we have the following 
results:

\begin{lemma}
\begin{enumerate}
\item The subspace $V^{p,n,\gamma}$ is spectrally equivalent
to $V^{p,n,(ij) \cdot \gamma}$ for any $i \neq j$. %\\
\item The subspace $V^{p,n,\gamma}$ is spectrally equivalent to
$V^{p,n,\sigma \cdot \gamma}$ for any $\sigma \in S_n$.
\end{enumerate}
\end{lemma}

The proof is straightforward, in light of the above discussion.

Before continuing, we note several basic facts about these
symmetry operators:
\begin{eqnarray}
& \chi_{ij}^2 & = {\rm Id} \label{E:Xsq} \\
& \chi_{ij} \chi_{ik} \chi_{ij} & = \chi_{jk} \label{E:Xtrans} \\
& \chi_{jk} U_{1j} & = U_{1k} \chi_{jk} \label{E:X&U} \\
& & \hbox{for any} \ i,j,k \ \hbox{such that} \ 2 \le i < j < k \le n
 \nonumber
\end{eqnarray}

From \eqref{E:Xsq}, we deduce that $\chi_{23}$ has only +1 and -1 as eigenvalues.
From \eqref{E:Xtrans}, we see that given $\chi_{23}, \chi_{24}, \ldots, \chi_{2n}$,
we can generate (by composition) any other $\chi_{ij}$ (for $2 \le i < j \le n$).
Equation \eqref{E:X&U} will be useful in the next section.

\section{Subspaces and sub-subspaces}

We are now able to use Lemma \ref{L:Ucomm} (\ref{L:gamma_beta}) to divide up the
eigenvalues of $\triangle_{p,n}(k)$ on $V^{p,\gamma}$ for any $\gamma$,
using Theorem \ref{T:red}, to be proved shortly. However, we first need to define
certain maps for convenience.

\begin{definition}
We define a map from ${\mathbb Z}^n$ to
${\mathbb Z}^{n-1}$ which omits the $i$th index and is denoted $p_i$:
$$p_i: (a_1,\ldots,a_{i-1},a_i,a_{i+1},\ldots,a_n) \mapsto
(a_1,\ldots,a_{i-1},a_{i+1},\ldots,a_n).$$
This then induces a projection $p_i^*$ from $\mathcal{F}^{-k}_n$
onto $\mathcal{F}^{-k}_{n-1}$, defined to be the linear operator with
action on the basis elements given by:
$$p_i^*: \psi_{\beta}(k) \mapsto \psi_{p_i(\beta)}(k).$$
\end{definition}
We can now state the theorem.

\begin{theorem} \label{T:red}
\begin{enumerate}
\item For any
multi-index $\gamma$ such that
$\gamma_n \ge 2$, $V^{p,n,\gamma}$ is spectrally equivalent to
$$(\ker U_{1n} \cap V^{p,n,\gamma'}) \oplus (\ker U_{1n} \cap
V^{p,n,\gamma''}) \oplus V^{p,n,\gamma'''}$$
where $\gamma'=\gamma+(\gamma_n-1)e_1-(\gamma_n-1)e_n$ (so that
$(\gamma')_n = 1$) and $\gamma''=\gamma'+e_1-e_n$, $\gamma'''
=\gamma'+2e_1-2e_n$. %\\
\item \label{T:phi1} The subspace $V^{p,n,\gamma'''}$ is spectrally
equivalent to
$V^{p-1,n-1,p_n(\gamma''')}$. %\\series(current,x,m)
\item \label{T:phi2} The subspace $\ker U_{1n} \cap V^{p,n,\gamma'}$
is spectrally equivalent to %\\
$V^{p-1,n-1,p_n(\gamma')+e_1}$.
\end{enumerate}
\end{theorem}

\begin{proof}
         (i) By Lemma \ref{L:Ucomm} (\ref{L:gamma_beta}), the subspaces
         $V^{p,n,\gamma}$ and
         $V^{p,n,\gamma'}$ are spectrally equivalent.

         Now as mentioned before,
         $$V^{p,n,\gamma'} \cong {\im}\, U_{n1} \oplus \ker U_{1n}$$
         where each subspace is $\triangle(k)$-invariant (a consequence of
         Lemma \ref{L:Ucomm} (\ref{L:Ud})).
         But $\gamma_1 \ge 2$ by assumption; since we are adding a
         non-negative number to the first index of $\gamma$, we also have
         $(\gamma')_1 \ge 2$, which means that $U_{n1}$ (here a map from
         $V^{p,n,\gamma''}$ to $\im U_{n1}$) is one-to-one, and so
         $\im \, U_{n1}$ is spectrally equivalent to $V^{p,n,\gamma''}$.
         The above decomposition then implies that
         $$V^{p,n,\gamma} \ \hbox{is spectrally equivalent to} \,
         (\ker U_{1n} \cap V^{p,n,\gamma'}) \oplus
         V^{p,n,\gamma''}.$$
         Similarly we can decompose $V^{p,n,\gamma''}$ into
         $\im \, U_{n1} \oplus \ker U_{1n}$; again, $U_{n1}$ is one-to-one
         and thus an isomorphism from $V^{p,n,\gamma'''}$ to
         $\im \, U_{n1}$. Hence
         $V^{p,n,\gamma''}$ is spectrally equivalent to
         $(\ker U_{1n} \cap V^{p,n,\gamma''}) \oplus
         V^{p,n,\gamma'''}$.
         These two spectral equivalences then imply part (i) of the
         theorem.

         (ii) We are considering $V^{p,n,\gamma'''}$, where
         $(\gamma''')_n = -1$. Recall that $\psi_{\beta}$ is
         only defined if $\beta_i \ge 0$ for all $i$ from 1 to $n$,
         so any element of $V^{p,n,\gamma'''}$ must be of the form
         $v \wedge \tau^n$, for some $v$ in $V^{p-1,n,\gamma'''+e_n}$
         such that $i(Z_n)v=0=i(Z_{\bar{n}})v$.

         We construct a homomorphism
         \[ \varphi_1: V^{p,n,\gamma'''}
         \rightarrow V^{p-1,n-1,p_n(\gamma''')} \]
         which takes $v \wedge \tau^n$ to $p_n^*(v)$, where we extend
         $p_n^*$ by tensoring with the projection from $\Lambda^{p-1}
         ((\mathfrak{h}^{2n+1})^*)$ onto $\Lambda^{p-1}
         ((\mathfrak{h}^{2n-1})^*)$. It can easily be seen that this
         (linear) homomorphism $\varphi_1$ is in fact one-to-one and onto.

For $\theta_{j}$ and $\theta_{j}^{*}$ the operators defined in
Definition \ref{D:theta}, we have
         $\theta_n (v \wedge \tau^n) =0 = \theta_n^* (v \wedge \tau^n)$;
         further, for any $i,j = 1,2, \dots , n-1$, the operators
         $\theta_j \theta_i^*$ and $\theta_i^* \theta_j$ both commute
         with the homomorphism $\varphi_1$ (that is, $\varphi_1$ doesn't
         affect the action of these operators).
         Thus $\triangle(k)$
         commutes with $\varphi_1$, which proves part (ii).

To prove part (iii),
construct a linear mapping $\varphi_2$ from $\ker U_{1n} \cap V^{p,n,\gamma'}$
to
$V^{p-1,n-1,p_n(\gamma')+e_1}$ as follows.

Begin by specifying that $\varphi_2$ commutes with $e(\tau^j) a_j^*$ and $e(\tau^{\bj})
a_j$ for $2 \le j \le n-1$, and with $e(\tau^w)$, and also with
the adjoints of these operators. These operators map
$\ker U_{1n} \cap V^{p,n,\gamma'}$ to $\ker U_{1n} \cap V^{p+1,n,\gamma'}$
and $V^{p-1,n-1,p_n(\gamma')+e_1}$ to
$V^{p,n-1,p_n(\gamma')+e_1}$.

Proceed by defining a $2$-form, to be denoted $\omega_2$, in
$\ker U_{1n} \cap V^{p,n,\gamma'}$, by
\begin{equation*}
\omega_2 := \psi_{\gamma'-e_1-e_n} \tau^{\bar{1}} \wedge \tau^{\bar{n}}
\end{equation*}
and specify that $\varphi_2$ maps $\omega_2$ to $\psi_{p_n(\gamma')} \tau^{\bar{1}}$.

Define next forms $\omega_0, \omega_1, \omega_3$ in $\ker U_{1n} \cap
V^{*,n,\gamma',\alpha}$ to be the $1$-form, $2$-form and $3$-form respectively
given by
\begin{eqnarray*}
\omega_0 & := & (g+1)^{-1/2} \bigl( i(Z_{\bar{1}})a_1^* + i(Z_{\bar{n}})a_n^*
\bigr) \omega_2,
\\
\omega_3 & := & (g+1)^{-1/2} \bigl( e(\tau^1)a_1^* + e(\tau^n)a_n^*
\bigr) \omega_2 \hspace{5mm} \text{and}  \\
\omega_1 & := & (g+2)^{-1/2} \bigl( e(\tau^1)a_1^* + e(\tau^n)a_n^*
\bigr) \omega_0
\end{eqnarray*}
for $g:=\gamma_1 + \gamma_n -1$;
these forms $\omega_j$ are still in $\ker U_{1n}$ since
$[e(\tau^1)a_1^*+e(\tau^n)a_n^*,U_{1n}]=0=[e(\tau^{\bar{1}})a_1+
e(\tau^{\bar{n}})a_n,U_{1n}]$.

Set $\varphi_2$ to also map $\omega_0$ to $\psi_{p_n(\gamma')+e_1}$,
map $\omega_3$ to $\psi_{p_n(\gamma')+e_1} \tau^1 \wedge \tau^{\bar{1}}$ and
map $\omega_1$ to $\psi_{p_n(\gamma')+2e_1} \tau^1$.
The mapping $\varphi_2$ can now be seen to be an isomorphism (for
example, by counting dimensions of the respective subspaces).

It's necessary to check that $\omega_0, \omega_1$ and $\omega_3$
all have length $1$ with respect to the inner product that we've chosen on
$V^{p,n,\gamma'}$, and also that $\omega_1$ is equal to
$-(g+2)^{-1/2} (i(Z_{\bar{1}}) a_1^* + i(Z_{\bar{n}}) a_n^*)
\omega_3$.
By considering the actions of adjoints of the operators discussed,
it follows (after some calculations) that
\begin{eqnarray*}
\varphi_{2} (\theta_{1}(k) + \theta_{n}(k)) = \theta_{1}(k) \varphi_{2}, \\
\varphi_{2} (\theta_{1}^{*}(k) + \theta_{n}^{*}(k)) = \theta_{1}^{*}(k)
\varphi_{2}
\end{eqnarray*}
on $\ker U_{1n} \cap V^{p,n,\gamma'}$, and that $\varphi_{2}$
commutes with the other operators involved in $d$ and $d^{*}$.

It follows that
$\varphi_2 \triangle_{p,n}(k) = \triangle_{p-1,n-1}(k) \varphi_2$ on all of
$\ker U_{1n} \cap V^{p,n,\gamma'}$.
Since we've already proved that $\varphi_2$ is an isomorphism, this completes
the proof.

 \end{proof}

      This theorem would also be true if we replaced $U_{1n}$ by
      $U_{12}$ or $U_{13}$
      and so on, so that instead of considering $\ker U_{1n} \cap V^{p,
      \gamma''}$, we need only consider the subspace $\ker U_{12} \cap
      \ldots \cap \ker U_{1n} \cap V^{p,(| \gamma |,0,\ldots ,0)}$
      (since all the eigenvalues ``missed'' here are eigenvalues of
      $\triangle_{p-1,n-1} (k)$).

      This subspace $\ker U_{12} \cap
      \ldots \cap \ker U_{1n} \cap V^{p,n, (| \gamma |,0,\ldots ,0)}$
      will be referred to as the reduced subspace,
      and denoted by $V^{p, n, |\gamma|}_{red}$.

      The subspace $V^{p,n,(| \gamma |,0,\ldots ,0)}$ is just
      $V^{p,n,\abs{\gamma}e_{1}}$, which is how it will be referred to
      from now on.

We now derive a basis for certain symmetric subspaces
and study the action of the Laplacian thereon.
As implied previously, the reduced subspace $V^{p,|\gamma|}_{red}$
can be further decomposed, this time with respect to the action of the
$\chi_{2j}$'s.

\begin{definition}
We call an element of this subspace which is also in the $+1$-eigenspace
of {\bf all} the $\chi_{2j}$'s (and thus of all $\chi_{ij}$'s,
by \eqref{E:Xtrans}) a {\bf symmetric} element, and an element in the $-1$-eigenspace
of {\bf any} $\chi_{ij}$ an {\bf anti-symmetric} element. These two
possibilities account for all of the subspace, i.e.
$$V^{p, |\gamma|}_{red} = (E_1 \chi_{23} \cap \ldots \cap
E_1 \chi_{2n}) \oplus (E_{-1} \chi_{23} + \ldots + E_{-1} \chi_{n-1,
n})$$
(where $E_{\lambda}A$
refers to the eigenspace of $A$ corresponding to the eigenvalue $\lambda$).
The symmetric subspace is defined to be $E_1 \chi_{23} \cap \ldots
\cap E_1 \chi_{2n} \cap \ker U_{12} \cap \ldots \cap \ker U_{1n}
\cap V^{p, \abs{\gamma}e_{1}}$ and will be denoted by
$V^{p, |\gamma|}_{symm}.$
\end{definition}

In order to understand the eigenvalues of the Laplacian on these
subspaces, we begin by
characterising explicitly all elements of  $E_{-1} \chi_{ij} \cap
V^{p,|\gamma|}_{red}$ for $2 \le i<j \le n$.
Note that this subspace is preserved by $\triangle_{p,n} (k)$.

To investigate the eigenvalues of the Laplacian on the anti-symmetric
subspace, we first prove a slightly more general lemma than is strictly
necessary.

\begin{lemma} \label{L:down_2}
For any multi-index $\hat{\gamma}$ such that $(\hat{\gamma})_{n-1}=0=
(\hat{\gamma})_n$, the subspaces $E_{-1} \chi_{n-1,n} \cap \ker U_{1,n-1}
\cap \ker U_{1,n} \cap V^{p,n,\hat{\gamma}}$ and $V^{p-2,n-2,p_{n-1}p_n
(\hat{\gamma})}$ are spectrally equivalent.
\end{lemma}

\begin{proof}
We construct a linear mapping $\varphi_3: E_{-1} \chi_{n-1,n} \cap
\ker U_{1,n-1} \cap \ker U_{1,n} \cap V^{p,n,\hat{\gamma}} \rightarrow
V^{p-2,n-2,p_{n-1}p_n (\hat{\gamma})}$, and show that $\varphi_3$ is
an isomorphism and commutes with the Laplacian, in a similar manner to
the proof of Theorem \ref{T:red}\ref{T:phi2}.

      Define a $3$-form, $\hat{\omega}_2$, by
      $$\hat{\omega}_2 := 2^{-1/2} \psi_{\hat{\gamma}-e_1}
      \tau^{\bar{1}} \wedge (\tau^n \wedge \tau^{\bar{n}} -
      \tau^{n-1} \wedge \tau^{\overline{n-1}});$$
      note that $\hat{\omega}_2 \in E_{-1} \chi_{n-1,n} \cap
      \ker U_{1,n-1} \cap \ker U_{1,n} \cap V^{p,n,\hat{\gamma}}$.
      Set $\varphi_3(\hat{\omega}_2) = \psi_{p_{n-1}p_n (\hat{\gamma})}$.

      Again as in the proof of Theorem \ref{T:red}\ref{T:phi2}, define the $2$-form
      $\hat{\omega}_0$, $3$-form $\hat{\omega}_1$, and $4$-form
      $\hat{\omega}_3$ as follows, for $\hat{g}:=(\hat{\gamma})_1$:
      \begin{eqnarray*}
      \hat{\omega}_0 & := & (\sqrt{\hat{g}+1})^{-1/2} \bigl(
      i(Z_{\bar{1}})a_1^* + i(Z_{\overline{n-1}})a_{n-1}^* + i(Z_{\bar{n}})a_n^*
        \bigr) \hat{\omega}_2,\\
\hat{\omega}_3 & := & (\sqrt{\hat{g}+1})^{-1/2} \bigl( e(\tau^1)a_1^* +
e(\tau^{n-1})a_{n-1}^* + e(\tau^n)a_n^*
\bigr) \hat{\omega}_2 \hspace{5mm} \text{and}  \\
\hat{\omega}_1 & := & (\sqrt{\hat{g}+2})^{-1/2} \bigl( e(\tau^1)a_1^* +
e(\tau^{n-1})a_{n-1}^* + e(\tau^n)a_n^* \bigr) \hat{\omega}_0
\end{eqnarray*}
      Note that the operators $e(\tau^1)a_1^* + e(\tau^{n-1})a_{n-1}^*
      + e(\tau^n)a_n^*$ and $e(\tau^{\bar{1}})a_1 + e(\tau^{\overline{n-1}})
      a_{n-1} + e(\tau^{\bar{n}})a_n$ both commute with $U_{1,n-1},
      U_{1,n}$ and $\chi_{n-1,n}$, as do their adjoints, so that
      $\hat{\omega}_0, \hat{\omega}_1$ and $\hat{\omega}_3$ are all in
      $E_{-1} \chi_{n-1,n} \cap \ker U_{1,n-1} \cap \ker U_{1,n}$.

      We set $\varphi_3$ to map $\hat{\omega}_0$ to
      $\psi_{p_{n-1}p_n(\gamma)+e_1}$,
$\hat{\omega}_3$ to $\psi_{p_{n-1}p_n(\hat{\gamma})+e_1}
\tau^1 \wedge \tau^{\bar{1}}$ and
$\hat{\omega}_1$ to $\psi_{p_{n-1}p_n(\hat{\gamma})+2e_1} \tau^1$.

      Again we must check that $\hat{\omega}_0, \hat{\omega}_1$ and
      $\hat{\omega}_3$
all have length $1$ with respect to the inner product that we've chosen on
$V^{p,n,\hat{\gamma},\alpha}$, and also that $\hat{\omega}_1$ is equal to
$-(g+2)^{-1/2} (i(Z_{\bar{1}}) a_1^* + i(Z_{\overline{n-1}})a_{n-1}^*
+ i(Z_{\bar{n}}) a_n^*)
\hat{\omega}_3$. If we then set $\varphi_3$ to commute with the operators
$e(\tau^j)a_j^*,i(Z_j)a_j,e(\tau^{\bar{j}})a_j,i(Z_{\bar{j}})a_j^*$
for $1 < j < n-1$ and $e(\tau^w), i(W)$, it follows as in the proof of
Theorem \ref{T:red} (\ref{T:phi2}) that $\varphi_3$ is an isomorphism
and commutes with $\triangle(k)$.
\end{proof}

\begin{corollary}
The subspace $E_{-1} \chi_{ij} \cap V^{p,n,\abs{\gamma}}_{red}$ is
spectrally equivalent to a subspace of $V^{p-2,n-2,(\abs{\gamma}+1)e_{1}}$
for $i \neq j, i,j=2,\dots,n$.
\end{corollary}

This follows from Lemma \ref{L:down_2} and appropriate use of the
symmetry operators.

Thus the eigenvalues of the Laplacian in the representation
$\bar{\beta}_{k}$ on the anti-symmetric subspace have already been
counted, in a sense, since they occur as eigenvalues of the Laplacian
on lower-degree forms on a lower-dimensional Heisenberg group.
Specifically, they are the same eigenvalues that we would get from applying
(say) the homomorphism $\varphi_{1}$ from Theorem \ref{T:red} (ii) twice.

We now begin to analyse the eigenvalues which occur on the symmetric subspace.

\begin{lemma} \label{L:symm_same}
The subspace $V^{p,n,|\gamma|}_{symm}$ is equal to
$E_1 \chi_{23} \cap \ldots \cap E_1 \chi_{2n} \cap \ker U_{12}
\cap V^{p,n,\abs{\gamma}e_{1}}$.
\end{lemma}

      The proof follows by repeated application of equation \eqref{E:X&U}.

We explicitly characterise all symmetric elements, beginning
this process by looking at $E_1 \chi_{23} \cap \ldots \cap E_1 \chi_{2n}
\cap V^{p,n,\abs{\gamma}e_{1}}$ (i.e. dropping the $\ker U_{12}$
condition).

\begin{definition}
We define a $p$-form, $\varepsilon (p,n)$, which is in $V^{p,n,\abs{\gamma}
e_1}$, by:
\begin{eqnarray*}
\varepsilon (p,n) & := & \left\{ \begin{array}{ll}
			\psi_{|\gamma| e_1} & \hbox{if $p$=0}
					\\
			\sum_{j=2}^n a_j^* e(\tau^j) \varepsilon (p-1,n) &
				\hbox{if $p$ is odd} \\
			(-2/p) \sum_{j=2}^n a_j e(\tau^{\bar{j}}) \varepsilon (p-1,n)
			& \hbox{if $p$ is even and $p \ge 2$}
			\end{array} \right.
\end{eqnarray*}
so that
\begin{eqnarray*}
& \varepsilon (1,n) = \sum_{j=2}^n \psi_{|\gamma|e_1+e_j}
\tau^j; \ \
& \varepsilon (2,n) = \sum_{j=2}^n \psi_{|\gamma|e_1} \tau^j
\wedge \tau^{\bar{j}}; \\
& \varepsilon (3,n) = \sum_{j,l=2}^n \psi_{|\gamma|e_1+e_j}
\tau^j \wedge \tau^l \wedge \tau^{\bar{l}}; \ \
& \varepsilon (4,n) = \sum_{j,l=2}^n \psi_{|\gamma|e_1}
\tau^j \wedge \tau^{\bar{j}} \wedge \tau^l \wedge \tau^{\bar{l}}
\ldots
\end{eqnarray*}
\end{definition}

Now $\varepsilon (p,n)$ is certainly in $E_1 \chi_{23} \cap \ldots \cap			
 E_1 \chi_{2n} \cap V^{p,n,\abs{\gamma}e_1}$.
 In fact, we have the following lemma.

\begin{lemma} \label{L:symmbas}
\begin{enumerate}
\item \label{L:alleps} For $p \ge 3$, a basis for $E_1 \chi_{23}
\cap \ldots \cap E_1 \chi_{2n} \cap V^{p,n,\abs{\gamma} e_1}$
is given by
\begin{eqnarray*}
& \left\{ \varepsilon(p,n), \; a_1^* \tau^1 \wedge \varepsilon(p-1,n), \; a_1
	\tau^{\bar{1}} \wedge \varepsilon (p-1,n), \, \tau^1 \wedge \tau^{\bar{1}}
	\wedge \varepsilon (p-2,n), \right. \\
& \tau^w \wedge \varepsilon (p-1,n), \; a_1^* \tau^w \wedge \tau^1 \wedge
	\varepsilon (p-2,n), \; a_1 \tau^w \wedge \tau^{\bar{1}} \wedge
	\varepsilon (p-2,n), \\
& \left. \tau^w \wedge \tau^1 \wedge \tau^{\bar{1}} \wedge \varepsilon (p-3,n)
	\right\}
\end{eqnarray*}
\item \label{L:redeps}
For $p$ even and $p \ge 4$, $p = 2q$ say, a basis of the symmetric
subspace is
\begin{eqnarray*}
& \left\{ -|\gamma| \tau^1 \wedge \tau^{\bar{1}} \wedge \varepsilon (2q-2,n)
+a_1 \tau^{\bar{1}} \wedge \varepsilon (2q-1,n), \right. \\
& \tau^1 \wedge \tau^{\bar{1}} \wedge \varepsilon (2q-2,n) +
\varepsilon (2q,n), \\
& \tau^w \wedge \tau^1 \wedge \tau^{\bar{1}} \wedge \varepsilon(2q-3,n)
+ a_1^* \tau^w \wedge \tau^1 \wedge \varepsilon(2q-2,n) + \tau^w
\wedge \varepsilon (2q-1,n), \\
& \left. a_1 \tau^w \wedge \tau^{\bar{1}} \wedge \varepsilon(2q-2,n)
\right\}
\end{eqnarray*}
\item The matrix of $\triangle_{2q,n}(k)$ acting on $V^{2q,|\gamma|}_{symm}$
with respect to the above basis
is then
\begin{eqnarray*}
& ( 2k\abs{\gamma}+k^2+nk)\id + \\
& \left( \begin{array}{cccc}
		(q-1)(n-q) & 0 & \sqrt{k} & \sqrt{k} \\
		-q \abs{\gamma} & q(n-q-1) & -q\sqrt{k} & 0 \\
		\sqrt{k} \abs{\gamma} & -\sqrt{k} & k+q(n-q) & 0 \\
		\sqrt{k}(\abs{\gamma}+n-q) & -\sqrt{k} & 0 & -k+q(n-q)
		\end{array}
		\right)
\end{eqnarray*}			
and its eigenvalues are
\begin{eqnarray} \label{E:sym_evals}
& \left\{ 2k\abs{\gamma}+nk+k^2+q(n-q), \right. \nn \\
& \left. 2k\abs{\gamma} + nk +k^2+ \frac{1}{2}n + q(n-q-1) 
\pm (\tfrac{1}{4}
n^2 + nk + 2k \abs{\gamma} +k^2)^{1/2} \right\}. 
\end{eqnarray}
where the first eigenvalue has multiplicity $2$.
\end{enumerate}
\end{lemma}

      \begin{proof}
      (i) It is easily checked that all these elements are in the
      required subspace. Also, it can be seen that the subspace
      $\ker i(Z_1) \cap \ker i(Z_{\bar{1}}) \cap \ker i(W) \cap E_1 \chi_{23}
\cap \ldots \cap E_1 \chi_{2n} \cap V^{p,n,\abs{\gamma}e_{1}}$
is spanned by $\varepsilon(p,n)$, i.e. $\varepsilon(p,n)$ is the
only symmetric element which doesn't contain $\tau^1, \tau^{\bar{1}}$,
or $\tau^w$. In this way we see that the given elements do in fact
span the subspace in question.

      (ii) The proof requires investigating the action of $U_{12}$
      (which is all that is necessary by Lemma \ref{L:symm_same})
      on linear combinations of the basis elements from (i).

      (iii) Proving this is a matter of (somewhat tedious) calculations,
      using firstly the formula \eqref{E:lap_rep} and secondly calculating
      the eigenvalues of the matrix, using say Maple.
      \end{proof}

We only need to find a basis for even $p$, due to the following
      fact:
a basis of the symmetric subspace for odd $p$ ($p=2q+1$ say) can be
derived from
the bases given in (\ref{L:redeps})
(more precisely, from the bases corresponding to $p=2q$
and $p=2q+2$) by judicious use of $e(\tau^w)$ and $i(W)$, since these
operators commute with $U_{ij}$ and $\chi_{ij}$ for all $i$ and $j$.

      Further, using the Hodge star operator (which could also
      give us the basis for $p$ odd from those for even $p$),
      we see that the
eigenvalues of $\triangle_{2q,n}(k)$ are exactly those of
$\triangle_{2n+1-2q,n}(k)$, so that we need not calculate the
      action of the Laplacian on odd forms separately.

We note specifically that all the eigenvalues
given in equation \eqref{E:sym_evals} are greater than $(n-1)k +
      k^2/c^2$.

      \begin{remark}
      Through tedious but straightforward calculations, one can
      verify that the eigenvalues of the Laplacian on symmetric $1$-forms
      and symmetric $2$-forms are also given by equation \eqref{E:sym_evals},
      though for $1$-forms the first eigenvalue has multiplicity
      one. Similarly, it can be shown that the eigenvalues
      of the Laplacian on $V^{p,n,0}_{symm}$ are just
      $ \{ k^2 + nk + q(n-q), k^2 + (n+1)k + q(n-q+1) \}$.
      For more details, see \cite{thesis}.
      \end{remark}

\section{The spectrum of the Laplacian in a representation}

In this section, we will prove the following result:

\begin{theorem} \label{T:all_evals}
For any positive integers $p, n$, let
$\triangle_{p,n}$
be the Laplacian on $p$-forms on $H^{2n+1}$ and $\triangle_{p,n}(k)$
the corresponding operator in the representation $\rho_k$ which corresponds
to Fourier transform in $k$ over the centre variable. Then the
eigenvalues of $\triangle_{p,n}(k)$
are
\begin{eqnarray*}
& \{ & 2k (g-1) + k^2 + (n-p)k, \\
% & & k^2 + (n-p+r)k + \left[\frac{r}{2}\right]\left(n-p+\left[
% \frac{r+1}{2}\right]\right), \\
& & k^2 + (n-p+r+1)k + \left\lfloor \frac{r+1}{2}\right\rfloor \left(n-p+
\left\lfloor \frac{r}{2} \right\rfloor +1\right), \\
& & 2k (g-1) + k^2 + (n-p+r)k +
 \left\lfloor \frac{r}{2}\right\rfloor \left(n-p+\left\lfloor 
 \frac{r+1}{2}\right\rfloor \right),\\
& & 2k g + k^2 + (n-p+r)k + \frac{1}{2}(n-p+r) + \left\lfloor \frac{r-1}
{2}\right\rfloor \left(n-p+\left\lfloor \frac{r}{2}\right\rfloor \right) \\
& & \pm (1/4(n-p+r)^2 + (n-p+r)k + 2k g + k^2)^{1/2}\\
& & : g \in {\mathbb Z}, g \ge 1 \ \hbox{and} \ r=1,\ldots,p \}.
\end{eqnarray*}
For any $k$, the lowest eigenvalue of $\triangle_{p,n}(k)$ is
$k^2 + (n-p)k$, and its
multiplicity is
$\binom{n}{p}$.
\end{theorem}

Here $\lfloor n \rfloor$ is defined to be the greatest integer smaller than
$n$.

The bulk of the work to prove Theorem \ref{T:all_evals} has already been done
(in particular, see equation \eqref{E:sym_evals}).
To complete the proof, note that the spectrum of $\triangle_{p,n}(k)$
contains all the eigenvalues of $\triangle_{p-1,n-1}(k)$. But this
latter set in turn includes all the eigenvalues of
$\triangle_{p-2,n-2}(k)$
and so on, so that the spectrum of $\triangle_{p,n}(k)$
contains the spectrum of $\triangle_{p-r,n-r}(k)$ for any $r$
between 1 and $p$.

The only new (as yet ``unlisted'') eigenvalues at each stage
(working from lower degree forms to higher degree)
are those which occur on the symmetric
subspaces $V^{p-r,n-r,\abs{\gamma}}_{symm}$, for $\abs{\gamma}=0,1,\dots$
. So these, together with the eigenvalues of the Laplacian
acting on functions, are all of the eigenvalues of the Laplacian
on $p$-forms, $\triangle_{p,n}(k)$.

It remains to prove that the given eigenvalue is indeed lowest, and that
its multiplicity is as specified.

\begin{lemma}
The lowest eigenvalue of $\triangle_{p,n}(k)$ is $(n-p)k + k^2$.
\end{lemma}

\begin{proof}
Note that this eigenvalue occurs; for example,
$\psi_{0}(k) \tau^1 \wedge \ldots \wedge \tau^p$
is an eigenvector with this eigenvalue.

      For most of the eigenvalues in the list given in Theorem \ref{T:all_evals},
      it is clear that they are greater than this eigenvalue (since
      we are considering the case $k > 0$).
      We need only consider the last eigenvalue which contains the
      negative square root. We have that
      $$2k g + k^2 + nk + \frac{1}{2}n
       - (\sqrt{1/4n^2 + nk + 2k g + k^2})
       > k^2 + (n-1)k,$$
      for $n,g \ge 1$. So this eigenvalue will always be greater
      than $k^2 + (n-p+r-1)k$, which is greater than or equal to
      $k^2 + (n-p)k$ since $r \ge 1$.

      So all other eigenvalues are greater than
      the eigenvalue in question.
\end{proof}

\begin{lemma} \label{L:low_mul}
The multiplicity of the lowest eigenvalue of
$\triangle_{p,n}(k)$ is $\tbinom{n}{p}$.
\end{lemma}

\begin{proof}
	We note from the proof of the preceding lemma
	that the eigenvalues of $\triangle_{p-r,n-p+r}(k)$ (for $r \le p$)
	on the symmetric subspace $V^{p-r,n-p+r,|\gamma|}_{symm}$
	are strictly greater than the value under consideration,
	$(n-p)k+k^2$, for any $|\gamma|$.
	
	This eigenvalue, which occurs as an eigenvalue of
	$\triangle_{0,n-p}(k)$ acting on $\psi_{0}(k)$ in
	$\mathcal{F}^{-k}_{n-p}$, is only found in
	the spectrum of $\triangle_{p,n}(k)$ due to repeated applications of
	Theorem \ref{T:red}\ref{T:phi1}, together with isometries $\chi_{j,n-p+r}$
(for certain values of $j$ between $1$ and $n-p+r$).
	
	That is, we have: $V^{0,n-p,0}$ is spectrally equivalent to
	$V^{1,n-p+1,\sigma_1 \cdot (-1,0,\ldots,0)}$
	(for some transposition $\sigma_1$ in $S_{n-p+1}$), which in turn is
	spectrally equivalent to $V^{2,n-p+2,\sigma_2 \cdot (-1,-1,0,
	\ldots,0)}$ (for $\sigma_2$ some permutation in $S_{n-p+2}$,
actually the product
	of $\sigma_1$ and a disjoint transposition) and so on; by induction,
	we infer
	that $V^{0,n-p,0}$ is spectrally equivalent to
	$V^{p,n,\sigma_p
	\cdot (-1,\ldots,-1,0,\ldots,0)}$ for some $\sigma_p \in S_n$
with the multi-index consisting of
	$-1$ repeated $p$ times and 0 repeated $n-p$ times.
	(Here $S_{q}$ stands for the permutation group on $q$ symbols.)

	There are
	$\binom{n}{p}$ different ways of choosing $\sigma_p$ which give
	different multi-indices, which proves this lemma, since
	all of these subspaces 
	$V^{p,n,\sigma_p \cdot -e_1-\ldots-e_p}$ are of
	(complex) dimension 1.
\end{proof}
	
	Note that the direct sum of these subspaces, the eigenspace of the
	lowest eigenvalue, is the subspace denoted by Lott in \cite{lott}
	by $\mathcal{S}^p$.
	
	This concludes the proof of the lemma and thus proves Theorem \ref{T:all_evals}.

\section{Calculation of Novikov-Shubin invariants}

In this section, we prove the following corollary.

\begin{corollary} \label{C:NSinv}
The $p$th Novikov-Shubin invariant of $H^{2n+1}$ is
given by
\begin{eqnarray*}
\alpha_p(H^{2n+1}) = \left\{ \begin{array}{ll}
							n+1, & p \ne n, n+1, \\
							\frac{1}{2}(n+1), & p = n, n+1.
							  \end{array}
					 \right.					
\end{eqnarray*}
Furthermore, for any discrete subgroup $\Gamma$ of $H^{2n+1}$ such that
$M=H/\Gamma$ is a compact manifold, $\alpha_p(M) = \alpha_p(H^{2n+1})$;
any manifold which is homotopy equivalent to such a manifold $M$ also
has the same Novikov-Shubin invariants, and any manifold $M'$ whose
fundamental group $\pi_1(M')$ is a discrete subgroup of $H^{2n+1}$
has its first Novikov-Shubin invariant given by
$\alpha_1(M')=\alpha_1(H^{2n+1})$.
\end{corollary}

\begin{proof}
 This corollary follows from Theorem \ref{T:all_evals} by analysis of
      the form of ${\rm Tr}_{\Gamma} e^{-t \triangle_p}$ and by
      showing that the lowest eigenvalue of $\triangle_p(k)$ does
      indeed determine the decay, as follows.

We need to first establish how ${\rm Tr}_\Gamma e^{-T\triangle_p}$ can be
found given only the eigenvalues of $\triangle_p(k)$ for all $k$.

 Let $\Gamma$ be a discrete subgroup of $H$, and let $\mathcal{A}$ be
      the von Neumann algebra defined as in section 1.
      Suppose we have an operator $A$ in $\mathcal{A}$ which is also
      $\Gamma$-trace class, is positive, self-adjoint, and has smooth
      kernel $k_A(x,y)$. Suppose also that $A$ is not only
      $\Gamma$-invariant, but also $H$-invariant.

      Then $L_g A u = A L_g u$ for any $g \in H$, $u \in L^2(H)$.
      This implies that the kernel of $A$ satisfies
      $k_A(g^{-1} x, y) = k_A(x,g y)$ for all $x,y \in H$,
      i.e. that $k_A(x,y) = k_A(e,x^{-1} y)$ and $k_A$ is a
      convolution kernel. 

     In particular, if $A=\pi_{R}(f)$ for some $f \in 
      C^{\infty}_{0}(H)$, then 
      $${\rm Tr}_{\Gamma}(A) = {\rm vol}(H/\Gamma) t(A),$$
      where $t$ is the trace defined in Theorem \ref{T:AbPlThm}.
      Since $\pi_{R}(C^{\infty}_{0}(H))$ is dense in $\mathcal{V}_{H}$,
      we have that ${\rm Tr}_{\Gamma}(A)={\rm vol}(H/\Gamma) t(A)$
      for any $A \in \mathcal{V}_{H}$.
      (This argument is taken from \cite{a&s}; see that article for 
      more details. The assumption there is that the Lie group in 
      question is semi-simple, but the same arguments hold for
      nilpotent Lie groups, for example $H$.)   

      The Plancherel Theorem for $H$ (see \cite{dixmiervn} or
      \cite{follandnew})
      implies that this trace also
      decomposes under equation \eqref{E:PL4H},
      $$t(A) = \int_{k \in \mathbb{R}} {\rm tr}_k (A(k)) \abs{k}^n dk,$$
      where $A = \int^\oplus_{k \in \mathbb{R}} A(k) \abs{k}^n dk$.
      We take $A(k)$ to be the operator $A$ in the representation
      $\bar{\beta}_k$ (since $A$ is left-invariant), so that
      ${\rm tr}_k$ is just the Hilbert-Schmidt trace on trace-class
      operators on $\bar{\mathcal{F}}^{k}_n$.

      In particular, if $\{ \lambda_j (k) \}_{j=1}^\infty$ are the eigenvalues
      of $A(k)$, with $\lambda_1 (k) \ge \lambda_2 (k) \ge \dots$, then
      ${\rm tr}_k (A(k)) = \sum_{j=1}^\infty \lambda_j (k)$.

      Now to apply this theory to our situation. We have to extend
      all the traces above by tensoring with the trace on
      $End(\Lambda^p(\mathfrak{h}^*))$, but this carries through all
      the above discussion.
      The Laplacian $\triangle_{p,n}$
      is (left) $H$-invariant, positive, self-adjoint,
      and has smooth kernel; therefore, so does the heat operator
      $e^{-t \triangle_{p,n}}$. Further, the heat operator is bounded,
      and thus is in $\mathcal{A}$.

      We still need the eigenvalues of the heat operator; however, it
      is a corollary of the spectral theorem for self-adjoint operators
      that if an operator $B$ has an eigenvalue $\lambda$, then
      $e^{-tB}$ will have an eigenvalue $e^{-t \lambda}$ (in fact, the
      eigenvector will be the same).

      So given all the eigenvalues of the Laplacian on $p$-forms on
      $H^{2n+1}$ in the representations
      $\bar{\beta}_k$, we can determine the Novikov-Shubin invariants.

      However, not all of this information is needed just to calculate
      the $p$th Novikov-Shubin invariant; only the value of the lowest
      eigenvalue of $\triangle_p(k)$ and its multiplicity (for all $k$)
      is strictly necessary, as shown below.

      Here we consider $k$ to be fixed.
      Suppose we have ordered the eigenvalues of $\triangle_p(k)$; that is,
      (indexing by positive integers) so that
      $\lambda_1(k) = \ldots = \lambda_m(k) < \lambda_{m+1}(k) \le
      \lambda_{m+2}(k) \ldots$
      (where $m$ is the multiplicity of the lowest eigenvalue).
      Then
      $$\sum_{j=1}^\infty e^{-T \lambda_j(k)} = m e^{-T \lambda_1(k)}
      + \sum_{j=m+1}^\infty e^{-T \lambda_j(k)}$$
      As $\triangle_{1,n}(k)$ is an unbounded operator, we have
      $\lambda_j(k) \rightarrow \infty$ as $j \rightarrow \infty$;
      thus the last term will always be less than
      $e^{-T \lambda_1(k)}$,
      and will not affect the decay as $T \rightarrow \infty$.

      A further complication occurs if the multiplicity of the lowest
      eigenvalue of $\triangle_{p,n}(k)$ varies with $k$; however, this is
      not the case here. As proved above, the multiplicity of the lowest
      eigenvalue of the Laplacian on $p$-forms on the Heisenberg group
      depends only on $n$ and $p$.

      Finally we need to consider the relevant integrals. For $a>0$ and
      $b>0$ constants, it can be shown
      %(see for example \cite{g&r})
      that
      $$\int_0^\infty k^m e^{-T(ak+k^2.f(k))} dk =
      \left( \frac{1}{aT} \right)^{m+1} + O(T^{-m-2}),$$
      $$\int_0^\infty k^m e^{-T(bk^2+k^3.f(k))} dk =
      O(T^{-1/2(m+1)})$$
      for $m$ a positive integer, and $f(k)$ a power series in $k$
      (positive for all $k>0$).

We see that the Heisenberg group $H^{2n+1}$ is $L^2$ acyclic, that is
has all $L^2$ Betti numbers zero, since $\lim_{T \rightarrow \infty} 
{\rm Tr}_{\Gamma} (e^{-T \triangle_{p,n}}) =0$. In fact, this was known
previously; it can be shown using techniques in \cite{c&m_acyclic}.

      So if $p<n$, we have that the decay of
      $${\rm Tr}_\Gamma e^{-T \triangle_p}
      = \int_{\mathbb R} {\rm tr}_k e^{-T \triangle_p(n,k)} |k|^n dk$$
	  is determined by the lowest eigenvalue, so that the integral is
	  of the first kind, with $m=n$ and $a=n-p$ (with $f(k)=1$).
	  Thus $\alpha_p(H^{2n+1}) = n+1$ if $p<n$, and using the Hodge
	  star operator, if $p>n+1$.
	
	  If $p=n$, the decay of the trace of the heat kernel is again
	  determined by the lowest eigenvalue $k^2$, and
	  thus by an integral of the second kind, with $m=n$ and $b=1$
	  (with $f(k)=0$). Thus $\alpha_n(H^{2n+1}) = 1/2(n+1) =
	  \alpha_{n+1}(H^{2n+1})$, again using the Hodge star operator.
	
	  The remainder of Corollary \ref{C:NSinv} follows from the definition of
	  the Novikov-Shubin invariants and the fact that they are
	  homotopy-invariant (see \cite{bmw, g&s}), as well as the fact that
	  the first Novikov-Shubin invariant $\alpha_1(M)$ of a manifold
	  $M$ a function only of the fundamental group of $M$,
which was proved in \cite{l&l}
	  \end{proof}

\section{General nilpotent Lie groups}

For all of the following sections, the universal reference for background
and definitions is
\cite{c&g}, which covers these topics in detail; this reference will
be assumed even if not specifically mentioned.

This section contains few new results; in particular, the results given
here agree
with the calculations in \cite{pesce} on the spectrum of the Laplacian on functions in a
representation, and are mostly an
elaboration of Appendix A of that article. However, the definitions of this
section are necessary for later sections.

We first outline the notation that will be used frequently from now on.

Let $\mathfrak{n}$ be a general nilpotent Lie algebra, with $N$ the
corresponding
connected and simply connected Lie group. %uniquely associated with it.
(Again, $N$ is unique up to isomorphism; see \cite{c&g}.)
Let $\mathfrak{z}$ be the centre of $\mathfrak{n}$, and $\mathfrak{v}$
the complement of $\mathfrak{z}$. Let $l$ be the dimension of
$\mathfrak{z}$, and $m$ the dimension of $\mathfrak{v}$.

\subsection*{Explicit formulae for $d$, $d^*$}
      From now on, we'll consider $N$ to be a step 2 nilpotent
      Lie group.

      Let $\{ X_{1}, \dots, X_{m+l} \}$ be an orthonormal basis for
$\mathfrak{n}$, %left transl etc.
and let $X_j$ also denote the left-invariant vector field derived
from $X_j$.
Let $\{ \tau^{1}, \dots, \tau^{m+l} \}$ be the corresponding
basis of 1-forms.
Then with respect to these bases, we can find explicit
formulae for $d, d^{*}$ and the Laplacian on functions and 1-forms at least.

If we select a basis $X_1,\dots,X_{m+l}$ for $\mathfrak{n}$, and
let the structure constants be $C_{ij}^k$, that is
$[X_i,X_j]=\sum_k C_{ij}^k X_k$, then we have that
      \begin{eqnarray*}
d & = & \sum_{i=1}^{m+l} e(\tau^{i}) X_{i} + \sum_{i,j=1}^{m}
\sum_{k=1}^{l} C_{ij}^{k} e(\tau^{i}) e(\tau^{j}) i(X_{k}), \\
d^{*} & = & -\sum_{i=1}^{m+l} i(X_{i}) X_{i} + \sum_{i,j=1}^{m}
\sum_{k=1}^{l} C_{ij}^{k} e(\tau^{k}) i(X_{j}) i(X_{i}), \\
\triangle_0 & = & -\sum_{i=1}^{m+l} X_i^2, \\
\triangle_{1} & = & -\left( \sum_{i,j=1}^{m+l} [X_{i}, X_{j}]
e(\tau^{i})i(X_{j}) + X_{j}^2 - \sum_{k} X_{j}
C^{k}_{ji} e(\tau^{i}) i(X_k) \right. \\
& & \left. + \sum_{k>j} C_{jk}^{i} X_{j} e(\tau^{i})i(X_{k})
-C_{jk}^{i} \sum_{q} C_{jk}^{m} e(\tau^{i})i(X_{m}) \right).
\end{eqnarray*}

      If instead we allow $\{ X_1, \dots, X_{m+l} \}$ to be complex
      vector fields (or to be an algebraic basis for $u(\mathfrak{n})$ ),
      which are orthonormal with respect to our chosen basis,
      then we have slightly different formulae for $d^*, \triangle_0$ and 
$\triangle_1$:
\begin{eqnarray}
d^{*} & = & -\sum_{i=1}^{m+l} i(X_{i}) \bar{X}_{i} + \sum_{i,j=1}^{m}
\sum_{k=1}^{l} \bar{C}_{ij}^{k} e(\tau^{k}) i(X_{j}) i(X_{i}), \nn \\
\triangle_0 & = & -\sum_{i=1}^{m+l} \bar{X}_i X_i, \label{E:lap0} \\
\triangle_{1} & = & -\left( \sum_{i,j=1}^{m+l} [X_{i}, \bar{X}_{j}]
e(\tau^{i})i(X_{j}) + \bar{X}_{j} X_{j} - \sum_{k} \bar{X}_{j}
C^{k}_{ji} e(\tau^{i}) i(X_k) \right. \nn \\
& & \left. + \sum_{k>j} \bar{C}_{jk}^{i} X_{j} e(\tau^{i})i(X_{k})
-\bar{C}_{jk}^{i} \sum_{q} C_{jk}^{m} e(\tau^{i})i(X_{m}) \right).
\label{E:lap_on_N}
\end{eqnarray}

      \subsection*{Kirillov theory}
Take any element $\lambda \in \mathfrak{n}^{*}$. Then we can define
a character $\zeta_\lambda$
on $Z$ (the centre of $N$, and the image of $\mathfrak{z}$
under $\exp$) to be
$\zeta_{\lambda} (\exp z) = e^{i \lambda(z)}$ for $z \in Z$.

Let $\pi_{\lambda}$ be the representation of $N$ induced
(in the sense of Mackey) from this representation $\zeta_\lambda$ of $Z$.
We write $\mathcal{H}_\lambda$ for its representation space.
We also denote the corresponding representation of $\mathfrak{n}$ by
$\pi_{\lambda}$.

      Kirillov theory (see for example \cite{c&g, kirillov}) tells us
      that every unitary representation $\pi$ of $N$ is unitarily
      equivalent to $\pi_\lambda$ for some $\lambda$; furthermore,
      two representations $\pi_\lambda, \pi_{\lambda'}$ are
      unitarily equivalent iff they are in the same ${\rm Ad}^* (N)$
      orbit, i.e. if there exists an element $g$ of $N$ such that
      $\lambda = {\rm Ad}^* g (\lambda')$. That is, $\hat{N}$ is
      the set of coadjoint orbits of $\mathfrak{n}^*$.

In particular, for any elements $W$ of $\mathfrak{z}$ and $\lambda$
of $\mathfrak{n}^{*}$,
we have that
$$\pi_{\lambda} (W) = \sqrt{-1} \, \lambda(W) \id$$
since $\pi_{\lambda}$ is a unitary representation.

      We'll actually consider the conjugate representation
      $\bar{\pi}_\lambda$, since results from the Heisenberg group
      (where $\bar{\beta}_k$ was the relevant representation)
      will then be more easily comparable;
      it also corresponds to left-invariant operators.
 Again, this representation
      is canonically isomorphic to $\pi_{-\lambda}$, and so we'll
      write its representation space as $\mathcal{H}_{-\lambda}$;
      Note that for $W \in \mathfrak{z}$,
      \[ \bar{\pi}_{\lambda} (W) = -\sqrt{-1} \, \lambda(W) \id. \]

      \subsection*{The Plancherel theorem for nilpotent Lie groups}
      We define the bilinear form $b_\lambda$ on $\mathfrak{n}$
      associated to any $\lambda \in \mathfrak{n}^*$ as follows:
      $$b_\lambda(X,Y):=\lambda([X,Y]).$$
      We also define the radical of this bilinear form, $r_\lambda$, as
      $r_\lambda:=\{ X \in \mathfrak{n} : b_\lambda(X,Y) = 0 \;
      \forall Y \in \mathfrak{n} \}$.
      Then $b_\lambda$ is non-degenerate on $\mathfrak{n}/r_\lambda$;
      from the theory of linear algebra, we know that this space is
      even-dimensional, of dimension $2n$ say.

      \newcommand{\Pf}{\rm Pf}

      Let $\{ X_1, \dots,
      X_{2n} \}$ be a basis for $\mathfrak{n}/r_\lambda$. Then the
      Pfaffian $\Pf(\lambda)$ is defined, up to sign, by
      $$\Pf(\lambda)^2 = \det B_\lambda,$$
      where $B_\lambda$ is the matrix with $(i,j)$th entry
      $b_\lambda(X_i,X_j)$.
      Once a choice of sign is made, $\Pf(\lambda)$ is a polynomial
      function of $\lambda$; specifically, a polynomial in
      $\lambda_1, \dots, \lambda_{m+l}$ (where $\lambda_i=\lambda(X_i)$)
      of degree $n$ (see \cite{satake}).

      Then it is well-known (see for example \cite{c&g}) that
      the Plancherel measure on $\pi_\lambda$ is Lebesgue measure
      on $\hat{N}$
      multiplied by the Pfaffian $\Pf(\lambda)$; that is,
      \begin{equation} \label{E:Pl_on_N}
      L^2(N) \isom \int_{\hat{N}} \mathcal{H}_\lambda \otimes
      \mathcal{H}_{-\lambda} \abs{\Pf(\lambda)} d \lambda.
      \end{equation}

      Again, we take the Laplacian on $p$-forms on $N$, $\triangle_p$,
      to be left-invariant, and write $\triangle_p(\lambda)$ for its
      decomposition in the representation $\bar{\pi}_\lambda$.
      That is, $\triangle_p(\lambda)$ is an operator on
      $\mathcal{H}_{-\lambda} \otimes \Lambda^p(\mathfrak{n}^*)$.

      \subsection*{Lower bound on spectrum}
      For general nilpotent groups, we can in fact use a similar method
      to that of Lemma \ref{L:lower_bound}
      to find a lower bound on the spectrum of
      the Laplacian in a representation for any nilpotent Lie group,
not just a step 2 nilpotent Lie group.

      \begin{theorem}
      For any $\lambda \in \mathfrak{n}^*$,
      $\triangle_p(\lambda) \ge \abs{(\lambda|_{\mathfrak{z}})}^2 \id.$
      \end{theorem}

      \begin{proof}
      For $N$ and $\mathfrak{n}$ as above, let $\{ X_1, \dots, X_m \}$
      be a basis for $\mathfrak{v}$ and $\{ W_1, \dots, W_l \}$ a basis
      for $\mathfrak{z}$, with $\{ \tau^{W_1}, \dots, \tau^{W_l} \}$
      the dual basis. Identify these elements with left-invariant vector
      fields and 1-forms as before.

      Define the operators $d_z:= \sum_{q=1}^l e(\tau^{W_q}) W_q$
      and $d_v:=d-d_z$, which both take $L^2$ $p$-forms on $N$ to $L^2$
      $(p+1)$-forms on $N$.
      Now $d_v$ can be written
      \[ d_v = ( \sum_{j=1}^m e(\tau^j) X_j ) + \sum_{i,j,k}
      C_{ij}^k e(\tau^i) e(\tau^j) i(X_k) \]
      (where the last term implicitly includes the case of $X_{k}=W_q$,
      i.e. that $X_{k}$ is in the centre), but importantly,
      there is no term $e(\tau^{W_q})$ in $d_v$ (for any $q=1,\dots,l$).

      This implies that $i(W_q) d_v + d_v i(W_q) = 0$, which means that
      $d_z^* d_v + d_v d_z^* = 0$, since $d_z^* = -\sum_{q=1}^l
      i(W_q)W_q$. Similarly $d_v^* d_z + d_z d_v^* = 0$.

      So
      \begin{eqnarray*}
      \triangle_p & = & d_v^* d_v + d_v d_v^* + d_z^* d_z + d_z^* d_z \\
                  & = & d_v^* d_v + d_v d_v^* - \sum_{q=1}^l W_q^2 \\
                  & \ge & -\sum_{q=1}^l W_q^2
      \end{eqnarray*}
      where the inequality follows since $d_v^* d_v + d_v d_v^*$ is a
      positive operator. But this means that
      $ \triangle_p(\lambda) \ge \sum_{q=1}^l \lambda(W_q)^2. $
      \end{proof}

\section{Heisenberg-type groups}

The main reference for this section is \cite{c&g}.

Let $\mathfrak{n}$ be a step 2 nilpotent Lie algebra with positive definite
inner product $\innprod{.}{.}$. Let $\mathfrak{z}$ be the centre of
$\mathfrak{n}$, and let $\mathfrak{v}$ be the complement of $\mathfrak{z}$
in $\mathfrak{n}$.
For each element $W \in \mathfrak{z}$, define a skew-symmetric linear
transformation $J(W)$ from $\mathfrak{v}$ to $\mathfrak{v}$ by:
\[ \innprod{J(W)X}{Y} = \innprod{W}{[X,Y]} \]
for all $X,Y \in \mathfrak{v}$.

\begin{definition}
A step 2 nilpotent Lie algebra $\mathfrak{n}$ with metric
$\innprod{.}{.}$ is of Heisenberg type (or H-type) if
$J(W)^{2} = - \abs{W}^{2} \id$ on $\mathfrak{v}$ for all $W \in
\mathfrak{z}$.
\end{definition}

We can then derive the following formula:
\begin{equation}
\innprod{J(W)X}{J(W')X} = \innprod{W}{W'} \abs{X}^{2}
          \label{E:WW'} 
%\innprod{J(W)X}{J(W)Y} & = & \innprod{X}{Y} \abs{W}^{2} \\
%\abs{J(W)X} & = & \abs{X} \abs{W} \\
%J(W) \circ J(W') + J(W') \circ J(W) & = & -2 \innprod{W}{W'} \id
\end{equation}
which is true for all $W,W'$ in $\mathfrak{z}$, and for all $X$ in
$\mathfrak{v}$; this and other formulae concerning 
$J(W)$ can be found in, for example, \cite{cdkr}.

      There is a connection between $J(W)$ and representations
      $\pi_\lambda$; to see it more clearly, we'll need the following
      notation.

      \begin{definition}
      Let $\{ W_1, \dots , W_l \}$ be an orthonormal basis for
      $\mathfrak{z}$, and $\{ \tau^{W_1}, \dots, \tau^{W_l} \}$
      the dual basis for $\mathfrak{z}^*$.

      For any element $W = \sum_{q=1}^l A_q W_q$ of $\mathfrak{z}$
      (with $A_q \in \mathbb{R}$),
      define $$\lambda_W := \sum_{q=1}^l A_q \tau^{W_q},$$
      the corresponding element of $\mathfrak{z}^*$.

      Similarly, for any element $\lambda = \sum_{q=1}^l B_q \tau^{W_q}$
      of $\mathfrak{z}^*$, define the corresponding element $W_\lambda$ of
      $\mathfrak{z}$ by
      $$W_\lambda := \sum_{q=1}^l B_q W_q.$$
      Trivially, $\lambda_{W_\lambda} = \lambda$ and
      $W_{\lambda_W} = W$.
      \end{definition}

      Now by definition,
      \begin{equation} \label{E:J_lambda}
      \innprod{J(W_\lambda)U}{V}=\lambda([U,V]),
      \end{equation}
      for all $\lambda \in \mathfrak{z}^*$, $U,V$ in $\mathfrak{v}$.
      Equivalently,
      $$\innprod{J(W)U}{V}=\lambda_W([U,V])$$
      for all $W \in \mathfrak{z}, U,V \in \mathfrak{v}$, and we use
      these two equations interchangeably.

Useful for our purposes will be the following lemma, which has a
straightforward proof, but is not (as far as I know) found in the literature.

\begin{lemma} \label{L:X_and_Y}
Let $\mathfrak{n}$ be any step 2 nilpotent Lie algebra with positive
definite inner product $\innprod{.}{.}$. 
If $\mathfrak{n}$ is H-type, then for any nonzero $\lambda \in
\mathfrak{n}^{*}$,
there is a basis $\{ X_{j \lambda}, Y_{j \lambda} \}_{j=1}^{n}$
of $\mathfrak{v}$ such that
$$\lambda([X_{j \lambda}, X_{k \lambda}]) = 0 = \lambda([Y_{j \lambda},
Y_{k \lambda}])$$ %number?
for $j,k=1, \dots ,n$, and
$$\lambda([X_{j \lambda}, Y_{k \lambda}]) = \delta_{jk}
\abs{\lambda}.$$
We take the inner product on $\mathfrak{n}^{*}$ to be that
induced by the inner product on $\mathfrak{n}$.
\end{lemma}

\begin{proof}
Choose any non-zero $\lambda \in
\mathfrak{n}^{*}$; in fact, we will assume without loss of generality
that $\lambda \in \mathfrak{z}^{*}$ (replacing $\lambda$ by another element
in its $\text{Ad}^{*}(N)$ orbit if necessary).

Now $\{ \mathfrak{v}, b_{\lambda} \}$ is a symplectic vector space
(because $b_\lambda$ is an anti-symmetric bilinear form on $\mathfrak{v}$, 
which is non-degenerate
since $\mathfrak{n}$ is H-type). So we can find a basis for $\mathfrak{v}$
(which depends on $\lambda$) $u_{1 \lambda}, \dots, u_{n \lambda}, v_{1 \lambda},
\dots, v_{n \lambda}$ such that 
$$b_\lambda (u_{i \lambda}, u_{j \lambda}) = 0 = b_\lambda (v_{i \lambda}, 
v_{j \lambda}) \ \mbox{and} \ 
b_\lambda (u_{i \lambda}, v_{j \lambda})=\delta_{ij}.$$
(For the proof, and more on symplectic vector spaces, see \cite{guill&stern}.) 
However, these elements $u_{i \lambda}, v_{j \lambda}$ are not necessarily
normalized. We define $X_{j \lambda}:=u_{j \lambda}/\norm{u_{j \lambda}}$,
and $Y_{j \lambda}:=v_{j \lambda}/\norm{v_{j \lambda}}$, so that
$X_{1 \lambda}, \dots , X_{n \lambda}, Y_{1 \lambda}, \dots, Y_{n \lambda}$
are an orthonormal basis for $\mathfrak{v}$.

Then since $\mathfrak{n}$ is H-type, we have that
\begin{eqnarray*}
\innprod{J(W_\lambda)^2 X_{i \lambda}}{X_{i \lambda}} & = & - \abs{\lambda}^2 \\
\implies -\innprod{J(W_\lambda)X_{i \lambda}}{J(W_\lambda)X_{i \lambda}} & = &
-\abs{\lambda}^2 \\
\implies J(W_\lambda)X_{i \lambda} & = & \abs{\lambda} Y_{i \lambda},
\end{eqnarray*}
where the last implication follows because $X_{i \lambda}$ is a scalar multiple
of $u_{i \lambda}$ (and because of the equation \eqref{E:J_lambda} which 
connects $b_\lambda$ and $J(W_\lambda)$).

\end{proof}

In fact, the converse of this lemma is also true; if such a basis 
of $\mathfrak{z}$ exists for any nonzero $\lambda \in \mathfrak{n}^*$,
then $\mathfrak{n}$ is H-type (see \cite{thesis}).

\begin{corollary} \label{C:X_Y_cr}
For any $\mathfrak{n}, \lambda$ as above, the basis $X_{j \lambda}, Y_{j
\lambda}$ of $\mathfrak{v}$ satisfies
$$[X_{j \lambda}, Y_{j \lambda}] = \frac{1}{\abs{\lambda}} \sum_{q=1}^l
\lambda_{q} W_{q}.$$
\end{corollary}

\begin{proof}
Define $W_\lambda$ as before.
As noted in the proof of the preceding lemma, we have that
$J(W_\lambda)X_{j \lambda}=\abs{\lambda} Y_{j \lambda}$.
But from equation \eqref{E:WW'}, we have that for any $p =1, \dots
,l$:
\begin{eqnarray*}
\innprod{J(W_{p})X_{j \lambda}}{J(W_\lambda)X_{j \lambda}} & = &
    \innprod{W_{p}}{W_\lambda} \\
\implies \innprod{J(W_{p})X_{j \lambda}}{\abs{\lambda} Y_{j \lambda}}
    & = & \lambda_{p} \\
\implies \abs{\lambda}\innprod{W_{p}}{[X_{j \lambda},Y_{j \lambda}]}
    & = & \lambda_{p}
\end{eqnarray*}
But this is true for all $p$, so the result follows.
\end{proof}

      \begin{remark}
This corollary says nothing about other commutation relations, such as
$[X_{j \lambda}, X_{k \lambda}]$; indeed, the only groups for which all
other commutation relations vanish are the Heisenberg groups.
\end{remark}

\begin{definition}
For any $\lambda \in \mathfrak{n}^* / \{0 \}$, we define 
$Z_{j \lambda}$ and $Z_{\bar{j} \lambda}$ to be the elements of
$u(\mathfrak{n})$ given by: 
\[ Z_{j \lambda}:=2^{-1/2} (X_{j \lambda} - i Y_{j \lambda}), \quad
Z_{\bar{j} \lambda}:=2^{-1/2} (X_{j \lambda} - i Y_{j \lambda}). \]
\end{definition}

The commutation relations of these elements are
\[ [Z_{j \lambda}, Z_{\bar{j} \lambda}] = i \abs{\lambda}^{-1} 
\sum_{q=1}^l \lambda_q W_q, \]
from Corollary \ref{C:X_Y_cr}. 
Thus $\bar{\pi}_\lambda ([Z_{j \lambda}, Z_{\bar{j} \lambda}]) = 
\abs{\lambda}$.

We can also think of $Z_{j \lambda}, Z_{\bar{j} \lambda}$ as complex 
left-invariant vector fields acting on $N$. With respect to them,
we can write
$$ \triangle_0(\lambda) = -\sum_{j=1}^n \left( Z_{j \lambda} Z_{\bar{j} \lambda}
+ Z_{\bar{j} \lambda} Z_{j \lambda} \right) -\sum_{q=1}^l W_q^2. $$
In particular,
\[
\, [\triangle_0(\lambda), \bpiz{j}] = 2 \abs{\lambda} \bpiz{j}, \quad
\, [\triangle_0(\lambda), \bpizb{j}] = -2 \abs{\lambda} \bpizb{j} 
\]
so that $\bpiz{j}, \bpizb{j}$ act as raising and lowering operators with
respect to $\triangle_0(\lambda)$.

      \subsection*{Creation and annihilation operators}
Creation and annihilation operators and a complete basis for
$\mathcal{H}_{\lambda}$ can now be defined, analogously to their
definition for $\mathcal{F}^{k}_{n}$.

For any $j$ between $1$ and $n$, let $a_{j}, a_{j}^{*}$ be the operators on
$\mathcal{H}_{-\lambda}$ defined by:
\[ a_{j} := \sqrt{-1} \abs{\lambda}^{-1/2} \bar{\pi}_{\lambda} (Z_{\bar{j}
\lambda}), \quad
 a_{j}^{*} := \sqrt{-1} \abs{\lambda}^{-1/2} \bar{\pi}_{\lambda} (Z_{j
\lambda}). \]
Then $[a_{j}, a_{j}^{*}]=\id$. We call $a_{j}$ an annhilation operator
and $a_{j}^{*}$ a creation operator.

For any $\lambda$, choose an element $v \in \mathcal{H}_{-\lambda}$
% actually in H^{\infty}_{\lambda}
which is in the kernel of $a_j$ for
all $j=1,\dots,n$.
(This element is unique up to scalar multiples,
otherwise the following construction would give a subspace of
$\mathcal{H}_{\lambda}$ which was $\bar{\pi}_{\lambda}$-invariant;
but this is impossible since $\bar{\pi}_{\lambda}$ is an irreducible
representation.) Define $\psi_{0}(\lambda)$ to be $v/\norm{v}$.

For any multi-index $\beta \in \mathbb{Z}_{+}^{n}$, we define
$$\psi_{\beta} (\lambda) := \frac{1}{\sqrt{\beta !}} (a^{*})^{\beta}
\psi_{0}(\lambda).$$

Then $\{ \psi_{\beta}(\lambda) \}_{\beta \in \mathbb{Z}^{n}_{+}}$ is
a complete basis of $\mathcal{H}_{\lambda}$ - otherwise, again, it
would be the basis for a closed, $\bar{\pi}_{\lambda}$-invariant subspace of
$\mathcal{H}_{\lambda}$.

      This leads to an explicit realisation of the representation
      $\bar{\pi}_\lambda$, with representation space $\mathcal{F}^{-\lambda}$,
      the generalized anti-Fock space (i.e. the conjugate to the
      generalized Fock space \cite{ricci}) - see \cite{thesis}.

%      This space is defined as
%      \begin{eqnarray*}
%      \mathcal{F}^{-\lambda} & := & \{ F: F \; \text{is defined and
%      anti-holomorphic on all of $\mathbb{C}^n$}, \\
%      & & \int_{\mathbb{C}^n} \abs{F(\bar{z})}^2 e^{-\abs{\lambda}
%      \abs{z}^2/4} dz < \infty \}.
%      \end{eqnarray*}

%      Then we set $a_j:=i \left( \tfrac{2}{\abs{\lambda}} \right)^{1/2}
%      \partial_{\bar{z}_j}$, and
%      $a_j^* := -i \left( \tfrac{\abs{\lambda}}{2} \right)^{1/2} \bar{z}_j$,
%      which gives
%      $$\psi_0(\lambda) = \left( \frac{\abs{\lambda}}{2 \pi}
%      \right)^{n/2}, \;
%      \psi_\beta(\lambda) = \left( \frac{\abs{\lambda}}{2 \pi} \right)^{n/2}
%      (-i)^{\abs{\beta}} \left( \frac{\abs{\lambda}}{2} \right)^{\abs{\beta}/2}
%      \frac{\bar{z}^\beta}{\sqrt{\beta !}}.$$

      \subsection*{An explicit formula for the Laplacian}
For H-type groups, the formulae for the Laplacian in particular
simplifies, so that we have

\begin{eqnarray} \label{E:H-type-lap}
\triangle_{1}(\lambda)
          & = & \abs{\lambda}^2 + n\abs{\lambda}+\sum_{j=1}^n \left( 2
           \abs{\lambda} a_j^*a_j +
          \abs{\lambda} (i(Z_j)e(\tau^j) +
             e(\tau^{\bj})i(Z_{\bj}) \right) \\
          & & + \sum_{i,j=1}^{n} \sum_{q=1}^{l}
          \left( \bpizb{j} C^{q}_{i,j} e(\tau^{i}) +
           \bpizb{j} C^{q}_{i+n,j} e(\tau^{\bar{i}}) \right. \nn \\
           & & \left. + \bpiz{j} C^{q}_{i,j+n} e(\tau^{i}) +
           \bpiz{j} C^{q}_{i+n,j+n} e(\tau^{\bar{i}}) \right)
           i(W_{q}) \nn \\
           & & + \left( \bpiz{j} \bar{C}^{q}_{i,j}i(Z_{i}) +
           \bpiz{j} \bar{C}^{q}_{i+n,j} i(Z_{\bar{i}}) \right. \nn \\
           & & \left. + \bpizb{j} C^{q}_{i,j+n} i(Z_{i}) +
           \bpizb{j} C^{q}_{i+n,j+n}
           i(Z_{\bar{i}}))e(\tau^{w_q}) \right)
             \nn 
\end{eqnarray}

There are similarities with the formula for the Laplacian on $p$-forms
on the Heisenberg group \eqref{E:lap_rep}, but the middle terms
(involving $C_{i,j}^{q}$ and so on) are rather different.
Nevertheless, we can list some of the eigenvalues of this
Laplacian in a representation, using these similarities.

\begin{lemma} \label{L:3evalues}
For any multi-index $\beta \in \mathbb{Z}_{+}^{n}$, we define
the following elements of $\mathcal{H}_{-\lambda} \otimes
\mathfrak{n}^{*}$:
\[
v_{1}:=\sum_{j=1}^{n} \sqrt{\beta_j+1} \, \psi_{\beta+e_j} (\lambda)
\tau^{j}, 
v_{2}:=\sum_{j=1}^{n} \sqrt{\beta_j}\, \psi_{\beta-e_j} (\lambda)
\tau^{\bar{j}}, 
v_{3}:=\sum_{q=1}^{l} \lambda_{p} \psi_{\beta} (\lambda)
\tau^{w_p}. 
\]
Then $\{ v_{1}, v_{2}, v_{3} \}$ span a
$\triangle_{1}(\lambda)$-invariant subspace of $\mathcal{H}_{-\lambda} \otimes
\mathfrak{n}^{*}$. Further, with respect to these elements,
$\triangle_{1}(\lambda)$ has matrix
\[
\triangle_{1}(\lambda) = (\abs{\lambda}(2\abs{\beta}+n) +
\abs{\lambda}^{2})\id + 
\left(
  \begin{array}{ccc}
  \abs{\lambda} & 0 & - \abs{\lambda}^{3/2} \\
  0 & -\abs{\lambda} & \abs{\lambda}^{3/2} \\
  -\abs{\lambda}^{-1/2} (\abs{\beta}+n) & \abs{\lambda}^{-1/2}
  \abs{\beta} & n
  \end{array}
\right)
\]
and eigenvalues
\[
\left\{ \abs{\lambda}(2\abs{\beta}+n) + \abs{\lambda}^2,  
\abs{\lambda}(2\abs{\beta}+n) + \abs{\lambda}^2 + \frac{n}{2}
\pm \sqrt{\frac{n^2}{4} + \abs{\lambda}(2\abs{\beta}+n) + \abs{\lambda}^2}
\right\}
\]
\end{lemma}
The proof is by computation, using the formula \ref{E:H-type-lap} for
$\triangle_1(\lambda)$.
      The matrix described in this lemma would be self-adjoint if the 1-forms
      $v_1, v_2, v_3$ were correctly normalized.
      Further, the first of the above eigenvalues comes from the action
      of $d$ on functions (i.e. the corresponding eigenvector is in
      $\im d$), but the other two do not.

      Note the similarities between this lemma and Lemma \ref{L:symmbas}
	(\ref{L:redeps})
      for $q=n$;
      in fact, if we set $c=1$, and identify $k$ with $\abs{\lambda}$
      and $\abs{\gamma}$ with $\abs{\beta}$, the eigenvalues agree
      exactly.

\subsection*{Symmetry operators on H-type groups}
In fact, H-type groups are easily classified. The following result was
noted by Kaplan in \cite{kaplan2}.

\begin{theorem} \label{P:Cliff}
The map $J:\mathfrak{z} \rightarrow {\rm End} (\mathfrak{v})$
extends to a representation of the Clifford algebra of $\mathfrak{z}$.
\end{theorem}

That is, $\mathfrak{v}$ is a Clifford module over $Cl(\mathfrak{z})$.

      Recall that the Clifford algebra associated to a vector space 
$V$ and quadratic form $q$, denoted $Cl(V,q)$, is generated by elements of $V$.
      For $v \in V$, we write $C(v)$ for Clifford multiplication by $v$
      (i.e. the corresponding element in $Cl(V,q)$). Then
      $C(v)C(w)+C(w)C(v)=-2q(v,w)$. For more on Clifford algebras, see
for example \cite{l&m}.

      An explicit example of the structure of H-type groups, as related to
      Clifford modules, is found in \cite{cdkr}, where the cases
      $\dim \mathfrak{z}=1, 3$ and $7$ are discussed. 

Now it is also well-known (see for example \cite{l&m}) that
for any finite-dimensional vector space $V$,
if $\dim V \cong 3 \pmod{4}$, then $Cl(V)$ has two non-isomorphic
irreducible representations, and that otherwise it only has one.
(All of these representations are finite-dimensional.)

We consider first the case that dim $\mathfrak{z}$ is not congruent to 3,
mod 4. Let $M$ be the unique (up to isomorphism) irreducible module
for $Cl(\mathfrak{z})$, and let $m$ be its dimension. Then
$$v \cong M_{1} \oplus \dots \oplus M_{r}$$
for some $r$, where each $M_{i}$ is a copy of $M$.
In particular, the action of $Cl(\mathfrak{z})$ is
the same on each $M_{i}$. That is, we can find a basis $\{ X_{j}
\}_{j=1}^{mr}$ of $\mathfrak{v}$ (where $X_{j} \in M_{i}$ iff
$m(i-1) < j \le mi$), such that
\begin{equation} \label{E:sameCRs}
[X_{mi+j}, X_{mi+l}]=[X_{j},X_{l}] 
\end{equation}
for all $j,l=1,\dots,m$,
and for all $i=1,\dots,r-1$.
The structure constants are similarly related: for all $i,j,l,q$
in the appropriate sets, we have that
$C^{q}_{mi+j,mi+l}=C^{q}_{j,l}$.

So whenever dim $\mathfrak{z}$ is not congruent to 3 mod 4, we have
the following definition.
\begin{definition}
The $(i,j)$ symmetry operator on $\mathfrak{n}$, $\chi_{ij}'$, is
defined for all $1 \le i < j \le r$ by the following rules:
$\chi_{ij}'$ is linear, $\chi_{ij}'$ maps $X_{m(i-1)+l}$ to
$X_{m(j-1)+l}$ and $X_{m(j-1)+l}$ to $X_{m(i-1)+l}$ for $l=1,\dots,m$,
and $\chi_{ij}'$ is the identity on the complement of $M_{i} \oplus
M_{j}$.
\end{definition}

(For example, if we think of $\mathfrak{v}$ as $M^{r}$,
then $\chi_{12}'(v_{1},v_{2},\dots,v_{r})$ would just be $(v_{2},
v_{1},\dots, v_{r})$.)

If $\dim \mathfrak{z} \cong 3 \pmod{4}$, then let $U$ and $V$ be
the non-isomorphic irreducible modules for $Cl(\mathfrak{z})$.
The complement $\mathfrak{v}$ must be isomorphic to
$U_{1} \oplus \dots \oplus U_{r} \oplus V_{r+1} \oplus \dots \oplus
V_{r+s}$, for some $r,s$, where $\mathfrak{z}$ acts on each $U_{i}$ as on
$U$ and on each $V_{j}$ as on $V$.

In this case, we define $\chi_{ij}'$ only for
$1 \le i < j \le r$ or for $r+1 \le i < j \le r+s$, but the rest of
the definition is the same.

\begin{lemma}
With notation as above, whenever $\chi_{ij}'$ is defined, it is a Lie
algebra isomorphism.
\end{lemma}

The proof is trivial, given \eqref{E:sameCRs}.

\section{The ``Double'' Heisenberg Group}
In this last section, we investigate a particular class of H-type
groups: those with 2-dimensional centre. We know from the above
classification that there
will be at most one such group of a given dimension (up to isomorphism);
in fact, since the irreducible modules of $Cl(\mathbb{R}^{2})$ are
4-dimensional, a H-type group with 2-dimensional centre must have
dimension $4n+2$, for some positive integer $n$.

\begin{definition}
Let $D^{4n+2}$ denote the ``double'' Heisenberg group of dimension
$4n+2$. That is, the Lie algebra $\mathfrak{d}^{4n+2}$ of $D^{4n+2}$
has basis $\{ X_{1}, \dots , X_{4n}, W_1, W_2 \}$ and non-zero commutation
relations defined by
\begin{eqnarray*}
\; [X_{4j+1},X_{4j+3}] = W_1, & & [X_{4j+1},X_{4j+4}] = W_2, \\
\; [X_{4j+2},X_{4j+3}] = W_2, & & [X_{4j+2},X_{4j+4}] = -W_1
\end{eqnarray*}
for $j=0,\dots,n-1$.
\end{definition}
(We also write $D$ instead of $D^{4n+2}$ and $\mathfrak{d}$ instead
of $\mathfrak{d}^{4n+2}$ when the dimension is understood.)

      \subsection*{Raising and lowering operators}
      We begin by defining the raising and lowering operators for
      $D^6$, then indicate how to generalise to $D^{4n+2}$.
      We fix a non-zero linear functional $\lambda \in \mathfrak{z}^*$
      throughout.

We find $Z_{1 \lambda}, \dots, Z_{\bar{2} \lambda}$ as indicated in
Lemma \ref{L:X_and_Y}; they are given by
      \begin{eqnarray*}
      Z_{1 \lambda} & := & (\sqrt{2} \abs{\lambda})^{-1}
                (i \lambda_1 X_1 + i \lambda_2 X_2 + \abs{\lambda} X_3), \\
      Z_{\bar{1} \lambda} & := & (\sqrt{2} \abs{\lambda})^{-1}
                (-i \lambda_1 X_1 - i \lambda_2 X_2 + \abs{\lambda} X_3), \\
      Z_{2 \lambda} & := & (\sqrt{2} \abs{\lambda})^{-1}
                (i \lambda_2 X_1 - i \lambda_1 X_2 + \abs{\lambda}
                X_4), \\
      Z_{\bar{2} \lambda} & := & (\sqrt{2} \abs{\lambda})^{-1}
                (-i \lambda_1 X_1 + i \lambda_2 X_2 + \abs{\lambda} X_4).
      \end{eqnarray*}
      These elements of $u(\mathfrak{d})$ have the following non-zero
      commutation relations:
      \begin{eqnarray*}
      \, [Z_{1 \lambda}, Z_{\bar{1} \lambda}] & = & i(\abs{\lambda})^{-1}
                (\lambda_1 W_1 + \lambda_2 W_2)
                = [Z_{2 \lambda}, Z_{\bar{2} \lambda}], \\
      \, [Z_{1 \lambda}, Z_{2 \lambda}] & = & i(\abs{\lambda})^{-1}
                (-\lambda_2 W_1 + \lambda_1 W_2)
                = -[Z_{\bar{1} \lambda}, Z_{\bar{2} \lambda}].
      \end{eqnarray*}
      But again, since $D$ is a H-type group, in the representation
      $\bar{\pi}_\lambda$ we have:
      \begin{eqnarray*}
      \bar{\pi}_\lambda([Z_{1 \lambda}, Z_{\bar{1} \lambda}]) & = \quad
                \abs{\lambda} & = \quad \bar{\pi}_{\lambda}
                ([Z_{2 \lambda}, Z_{\bar{2} \lambda}]), \\
      \bar{\pi}_{\lambda}([Z_{1 \lambda}, Z_{2 \lambda}]) & = \quad 0 &
                = \quad \bar{\pi}_{\lambda} ([Z_{\bar{1} \lambda}, Z_{\bar{2}
                \lambda}]).
      \end{eqnarray*}

      For $n \ge 2$, the remaining $Z_j$'s and $Z_{\bj}$'s are defined
      analogously; for example,
      $$Z_{3 \lambda} := (\sqrt{2} \abs{\lambda})^{-1}
                (i \lambda_1 X_5 + i \lambda_2 X_6 + \abs{\lambda} X_7).$$
      As we know, the commutation relations also carry over unchanged.

      We again write $\{ \tau^1, \tau^{\bar{1}}, \dots, \tau^{2n},
      \tau^{\overline{2n}} \}$ for the dual basis corresponding to \newline
      $\{ Z_{1 \lambda}, \dots, Z_{\overline{2n} \lambda} \}$.

      We define the creation and annihilation operators on
      $\mathcal{H}_{-\lambda}$ as we did for all H-type groups:
      $$ a_{j} = \sqrt{-1} \abs{\lambda}^{-1/2} \pi_{\lambda} (Z_{j
\lambda}), \quad
 a_{j}^{*} = \sqrt{-1} \abs{\lambda}^{-1/2} \pi_{\lambda} (Z_{\bar{j}
\lambda}).$$

      \subsection*{Commuting operators}
      For $D^{4n+2}$, the symmetry operators $\chi_{ij}'$ are easily
      defined; for example, $\chi_{12}'$ interchanges $X_1$ and $X_5$,
      $X_2$ and $X_6$, and so on, or equivalently, $Z_{1 \lambda}$ and
      $Z_{3 \lambda}$,
      $Z_{2 \lambda}$ and $Z_{4 \lambda}$, $Z_{\bar{1} \lambda}$
      and $Z_{\bar{3} \lambda}$, and
      $Z_{\bar{2} \lambda}$ and $Z_{\bar{4} \lambda}$ are interchanged.

      We can define transposition operators $U_{ij}$
      as for the Heisenberg group,
      in terms of $a_j$ and $e(\tau^j)$, but they do not commute with
      the Laplacian on 1-forms, $\triangle_1(\lambda)$.
      Instead, if $n \ge 2$,
      $[ \triangle_1(\lambda), U_{13} - U_{42}]=0 = [\triangle_1(\lambda),
      U_{31} - U_{24}];$
      and I conjecture that $\triangle_1(\lambda)$ also commutes with
      $U_{23}-U_{41}$ and thus with $U_{32}-U_{14}$.

      The Laplacian on 1-forms does not even commute with $U_{11}$ or
      $U_{22}$, but instead with $U_{11}-U_{22}$, so that there is
      no corresponding subspace $V^{1,n,\gamma}$, but instead two
      disjoint subspaces, as we'll see shortly.

      Use could be made of these operators in some way, but the
      situation is somewhat more complicated than for the Heisenberg group
      - due primarily to the extra non-zero commutation relations.

      \subsection*{Some eigenvalues of the Laplacian on 1-forms}
      For any multi-index $\beta \in \mathbb{Z}_+^n$ with all indices
      positive,
      define the following 1-forms for $j=1,\dots,n$:
      \begin{eqnarray*}
 u_j & := (a_{2j} e(\tau^{2j-1})
                - a_{2j-1} e(\tau^{2j})) \psi_\beta (\lambda), 
      v_j & := (a_{2j}^* e(\tau^{\overline{2j-1}})
                - a_{2j-1}^* e(\tau^{\overline{2j}})) \psi_\beta (\lambda), \\
      w_j & := (a_{2j-1}^* e(\tau^{2j-1})
                + a_{2j}^* e(\tau^{2j})) \psi_\beta(\lambda), 
      w_j' & := (a_{2j-1} e(\tau^{\overline{2j-1}})
                + a_{2j} e(\tau^{\overline{2j}})) \psi_\beta (\lambda).
      \end{eqnarray*}
Define also the number $\mu':=\abs{\lambda}(2n+2 \abs{\beta})+\abs{\lambda}^2$.

      \begin{theorem} \label{T:evalues_on_D}
      If $n \ge 2$, then for any multi-index $\beta \in \mathbb{Z}^n_+$,
      the Laplacian
      $\triangle_1 (\lambda)$ in the representation $\bar{\pi}_\lambda$
      acting on 1-forms on $D^{4n+2}$ has 
eigenvalues including $\{ \mu'-3 \abs{\lambda}, \mu'+3 \abs{\lambda}, \mu' + 
\abs{\lambda}, \mu'-\abs{\lambda} \}$, each with multiplicity $n-1$.
The corresponding eigenvectors are, respectively,
$\{ (\beta_{2j+1}+\beta_{2j+2})u_j-(\beta_{2j-1}+\beta_{2j})u_{j+1},
(\beta_{2j+1}+\beta_{2j+2})v_j-(\beta_{2j-1}+\beta_{2j})v_{j+1},
(\beta_{2j+1}+\beta_{2j+2})w_j-(\beta_{2j-1}+\beta_{2j})w_{j+1},
(\beta_{2j+1}+\beta_{2j+2})w_j'-(\beta_{2j-1}+\beta_{2j})w_{j+1}' \}$,
for $j=1, \dots, n-1$.
      There are also two $\triangle_1(\lambda)$-invariant
      subspaces with bases
$$
      \left\{ \sum_{j=1}^n u_j, \; \sum_{j=1}^n v_j, \ \lambda_2 \tau^{w_1}
                - \lambda_1 \tau^{w_2} \right\}, 
      \left\{ \sum_{j=1}^n w_j, \sum_{j=1}^n w_j', \; \lambda_1 \tau^{w_1}
                + \lambda_2 \tau^{w_2} \right\}.
$$
      With respect to these bases, $\triangle_1(\lambda)$ has matrices
      \begin{eqnarray*}
      & & (2\abs{\lambda}(n+\abs{\beta})+\abs{\lambda}^2) \id + 
      \left( \begin{array}{ccc}
                -3\abs{\lambda} & 0 & - \abs{\lambda}^{3/2} \\
                0 & 3\abs{\lambda} & \abs{\lambda}^{3/2} \\
                - \abs{\lambda}^{-1/2}\abs{\beta} &
                  -\abs{\lambda}^{-1/2} (\abs{\beta}+2n) & 2n
             \end{array}
      \right) \\ %\label{E:weird_matrix} \\
      & & \text{and} \\
      & & (2\abs{\lambda}(n+\abs{\beta})+\abs{\lambda}^2) \id +
      \left( \begin{array}{ccc}
                \abs{\lambda} & 0 & - \abs{\lambda}^{3/2} \\
                0 & -\abs{\lambda} & \abs{\lambda}^{3/2} \\
                - \abs{\lambda}^{-1/2}(\abs{\beta}+2n) &
                  -\abs{\lambda}^{-1/2} \abs{\beta} & 2n
             \end{array}
      \right).
      \end{eqnarray*}
      On the second subspace, $\triangle_1(\lambda)$ has eigenvalues
\[
      \left\{ 2\abs{\lambda}(n+\abs{\beta})+\abs{\lambda}^2, 
      2\abs{\lambda}(n+\abs{\beta})+\abs{\lambda}^2 +n \pm
           \sqrt{n^2+2\abs{\lambda}(n+\abs{\beta})+\abs{\lambda}^2)}
           \right\}.\]
      Define functions $\mu_{low}(b,n)$ and $\mu_{high}(b,n)$ for
      positive integers $b,n$ as follows:
\[
      \mu_{low}(b,n):=-\left( \frac{b+n+\sqrt{(b+n)^2+24n^2}}{2n} \right)
        \abs{\lambda}, \quad 
      \mu_{high}(b,n):=-3\abs{\lambda}
\]     
      Then the lowest eigenvalue $\mu_0$ of the first of the above
      matrices %\eqref{E:weird_matrix}
      is bounded by
      \[ \mu_{low}(\abs{\beta},n)+2(\abs{\beta} + n) \abs{\lambda}
      + \abs{\lambda}^2 < \mu_0 < \mu_{high}(\abs{\beta},n)+2 (\abs{\beta}
      + n) \abs{\lambda} + \abs{\lambda}^2,\]
      for $\abs{\lambda}>0$, $\abs{\beta},n \ge 1$,
      while the other two eigenvalues are greater than $2(\abs{\beta} + n)
      \abs{\lambda}+\abs{\lambda}^2$.
      \end{theorem}

      \begin{proof}
      Most of the proof consists of tedious calculations, using
      either of the
      formulae \eqref{E:lap_on_N} or \eqref{E:H-type-lap}.

      The eigenvalues of the first matrix in the theorem %\eqref{E:weird_matrix}
      are worth discussing in some detail, since they come from a
      cubic which is decidedly non-trivial to solve.

      Let $p(\mu)$ be the characteristic polynomial of this matrix
      (minus the constant term $2(\abs{\beta} + n)
      \abs{\lambda}+\abs{\lambda}^2$). That is,
      $$p(\mu) = \mu^{3}-2n\mu^{2} - \abs{\lambda}(2\abs{\beta} +
      9\abs{\lambda} + 2n)\mu + 12n \abs{\lambda}^{2}.$$
      Then we can approximate its
      zeros (i.e. the eigenvalues of the matrix) if we know where it
      is positive and negative. Calculations (for example, using a
      computer package such as Maple) give that
      $p(\mu_{low}(\abs{\beta},n))=\mu_{low}(\abs{\beta},n)^{3}-
      9 \mu_{low}(\abs{\beta},n) \abs{\lambda}^{2}$ is
      negative (since $\abs{\mu_{low}(\abs{\beta},n)}>3\abs{\lambda})$),
       while $p(\mu_{high}(\abs{\beta},n)) =
      6\abs{\beta} \abs{\lambda}^{2}$ is positive.
      Further, $p(0)$ is positive, while $p'(\mu)$ has a positive zero,
      indicating (by standard results in calculus) that the other two zeros
      of $p(\mu)$ are both positive.
      \end{proof}

      We briefly discuss special cases, i.e. what happens when some or
      all of the indices $\beta_i$ are zero.

      If $\beta_{2j-1}$, say, is zero, then $u_j$ and $w_j'$ both simplify;
      but if $\beta_{2j-1}=0=\beta_{2j}$, then $u_j$ and $w_j'$ are also
      zero.
      In particular, if $\beta=0$, then every $u_j$ and every $w_j'$
      are zero; also, $\abs{\beta}$ must be greater than or equal to 2 in
      order to have an eigenvector of the form
      $(\beta_{2l-1}+\beta_{2l}) u_j - (\beta_{2j-1}+\beta_{2j})u_l$
      (since one of $\beta_{2l-1}, \beta_{2l}$ must be non-zero,
      and one of $\beta_{2j-1}, \beta_{2j}$ must be non-zero).

      The case $\beta=0$ has to be considered separately, but it
      can be shown that in this case, all eigenvalues are greater than
      $2n \abs{\lambda} + \abs{\lambda}^2$.

      This motivates the following result.

      \begin{corollary} \label{C:D_lowest}
      The lowest eigenvalue of the Laplacian on 1-forms on
      $D^{4n+2}$ in the representation $\pi_\lambda$, $\triangle_{1,n}
      (\lambda)$, has multiplicity 1 for all $n,\lambda$, and lies between
      $\left( 2(n + 1) - \frac{n+1+\sqrt{(n+1)^2+24n^2}}{2n} \right)
      \abs{\lambda} + \abs{\lambda}^2$ and
      $( 2n - 1 ) \abs{\lambda} +
      \abs{\lambda}^2.$
      Further, the coefficient of $\abs{\lambda}$ in the lower bound
      is positive.
      \end{corollary}

      \begin{proof}
      For fixed $\beta$, the lowest eigenvalue on the first subspace
      in Theorem \ref{T:evalues_on_D} is between
      $$2 (\abs{\beta}+n) \abs{\lambda} + \abs{\lambda}^2 +
      \mu_{low}(\abs{\beta},n) \; \text{and} \;
      2 (\abs{\beta}+n) \abs{\lambda} + \abs{\lambda}^2 +
      \mu_{high}(\abs{\beta},n).$$
      Both $\mu_{low}(\abs{\beta},n)$ and $\mu_{high}(\abs{\beta},n)$
      are increasing as $\abs{\beta}$ increases; in particular,
      $ \mu_{low}(2,n) > \mu_{high}(1,n) \quad \forall n \ge 1.$

      All other eigenvalues are also greater than $\mu_{high}(1,n)+
      2(n+1)\abs{\lambda}+\abs{\lambda}^2$; in particular, the lowest
      of the other eigenvalues (coming from $(\beta_{2l-1}+\beta_{2l}) u_j
      - (\beta_{2j-1}+\beta_{2j})u_l$) is $(2n+1)\abs{\lambda}+
      \abs{\lambda}^2$, and the lowest eigenvalue on the second subspace
      is $2n\abs{\lambda}+\abs{\lambda}^2$.
      So the lowest eigenvalue of $\triangle_1(\lambda)$
      is between $\mu_{low}(1,n)+2(n+1)
      \abs{\lambda} + \abs{\lambda}^2$, and $\mu_{high}(1,n) + 2(n+1)
      \abs{\lambda} + \abs{\lambda}^{2}$.

      That the coefficient of $\abs{\lambda}$ is positive follows from
      more calculations. For $n=1$, the value of
      $\mu_{low}(1,1)+4 \abs{\lambda}$ is exactly
      $(3-\sqrt{7})\abs{\lambda}$;
      for $n>1$, we use the fact that $\sqrt{(n+1)^2+24n^2}$ is
      less than $5n+1$ to derive the estimate:
      $\mu_{low}(1,n)+(2n+2) \abs{\lambda}$ is greater than
      $(2n-1-\tfrac{1}{n}) \abs{\lambda}$, which is positive.
      \end{proof}

      \begin{corollary}
      For any $n \ge 1$, the first Novikov-Shubin invariant of
      $D^{4n+2}$ is given by
      $$\alpha_1(D^{4n+2}) = 2n+2 = \alpha_0(D^{4n+2}).$$
      \end{corollary}

      \begin{proof} From Corollary \ref{C:D_lowest}, we have an
      estimate for the lowest eigenvalue of
      $\triangle_{1,n}(\lambda)$, which has multiplicity of one for
      all $n$ and $\lambda$. As in section 4.4, we can now calculate
      the eigenvalues; most of the procedure of that section still
      holds here.  The result depends on the decay of the following
      integral:
      $$\int_{\mathbb{R}^2} e^{-T(a \abs{\lambda} + f(\abs{\lambda}).\abs{\lambda}^2)}
      \abs{\lambda}^{2n} d\lambda_1 d \lambda_2$$
      for $a$ positive and $f(x)$ a positive power series. We can rewrite this integral in polar
      coordinates; it becomes
      $$\int_{0}^{2 \pi} \int_0^\infty
      e^{-T(a r + f(r).r^2)} r^{2n+1} dr d \theta,$$
      which (again using an equation from ) evaluates to
      $2 \pi \left( \frac{1}{aT} \right)^{2n+2} + O(T^{-2n-3}).$
      \end{proof}
	
\appendix
	
\section{An explicit formula for the Laplacian}

The formula for the Laplacian in section 3.1 was given without proof -
though it was indicated how the explicit formulae for $d$ and $d^*$ could
be proved. Here we derive the formula for the Laplacian, given those
for $d$ and $d^*$.

First, we need to review some properties of the operators $e(*), i(*)$.

Let $U, V$ be %elements of $\mathfrak{h} \otimes \mathbb{C}$.
vectors selected from the basis $\{ Z_1, \ldots, Z_n,
Z_{\bar{1}}, \ldots, Z_{\bar{n}}, W \}$. Let $\tau^U, \tau^V$ be the
corresponding elements of the dual basis.
Then we have the following properties:
\begin{eqnarray}
\{ e(\tau^{U}), i(V) \} & = & \innprod{U}{V} \label{E:ACRfirst} \\
\{ e(\tau^{U}), e(\tau^{V}) \} & = & 0 = \{ i(U), i(V) \} \\
e(\tau^{V}) & = & [i(V)]^{*} \label{E:ACRlast}
\end{eqnarray}
where $\{.,.\}$ is the anti-commutator, $\{A,B\} := AB + BA$.
%explain tau^{V} etc.

We also note that vector fields such as $Z_j$,
which operate only on functions, commute with the operators
$e(\tau)$ and $i(V)$ for all $\tau, V$ in the above orthonormal bases.
%is this justified??

Finally, it can be shown that the adjoint of $Z_j$ is $-Z_{\bar{j}}$ and
the adjoint of $W$ is $-W$.

Recall from section 3.1 that
\begin{eqnarray*}
d & = & \sum_{j=1}^n \left( e(\tau^j)Z_j + e(\tau^{\bar{j}})Z_{\bar{j}} \right)
		+ e(\tau^w)W - i \sum_{j=1}^n e(\tau^j)e(\tau^{\bar{j}})i(W) \\
d^* & = & -\sum_{j=1}^n \left( i(Z_{\bar{j}})Z_j + i(Z_j)Z_{\bar{j}} \right) 
		- i(W) + i\sum_{j=1}^n e(\tau^w)i(Z_{\bar{j}})i(Z_j)
\end{eqnarray*}

We define, as before, the operators
$\theta_j$ for $j=1,\ldots,n$, by
\[ \theta_j = e(\tau^j)Z_j + e(\tau^{\bar{j}})Z_{\bar{j}}
                -ie(\tau^j)e(\tau^{\bar{j}})i(W) \]
so that $d=e(\tau^w)W + \sum_{j=1}^n \theta_j$. We can then define operators
$\eta_{j,l}$ and $A_j$ for $j \neq l$ and $j,l=1, \ldots ,n$:
\[ \eta_{j,l} := \theta_j \theta_l^* + \theta_l^* \theta_j, 
A_j		   := \theta_j \theta_j^* + \theta_j^* \theta_j;
\]
with these definitions, we can write
\[
\triangle_{p,n} = \sum_{j=1}^n A_j + \sum_{j \neq k} \eta_{j,k} -W^2.
\]
We now calculate $A_{j}$ and $\eta_{j,l}$.

Firstly,
\begin{eqnarray*}
A_{j} & = & \theta_{j} \theta_{j}^{*} + \theta_{j}^{*} \theta_{j} \\
      & = & \{ e(\tau^{j})Z_j + e(\tau^{\bar{j}})Z_{\bar{j}}
                -ie(\tau^j)e(\tau^{\bar{j}})i(W),
               -i(Z_{\bar{j}})Z_{j}-i(Z_{j})Z_{\bar{j}} \\
      & &          +ie(\tau^{w})i(Z_{\bar{j}})i(Z_{j}) \} \\
      & = & \{ e(\tau^{j})Z_{j},-i(Z_{j})Z_{\bar{j}} \} +
             ie(\tau^{w})i(Z_{\bar{j}})Z_{j} \\
      & & + \{ e(\tau^{\bar{j}}) Z_{\bar{j}}, -i(Z_{\bar{j}})Z_{j} \}
             -ie(\tau^{w})i(Z_{j})Z_{\bar{j}} \\
      & & -ie(\tau^{j})i(W)Z_{j}+ice(\tau^{\bar{j}})i(W)Z_{\bar{j}}
      \\
      & & + e(\tau^j)i(Z_j)e(\tau^{\bj})i(Z_{\bj})i(W)e(\tau^w) 
      + i(Z_j)e(\tau^j)i(Z_{\bj})e(\tau^{\bj})e(\tau^w)i(W)) \\
      & = & -2 Z_j Z_{\bar{j}} + iW \left( i(Z_j) e(\tau^j) + e(\tau^{\bar{j}})
      i(Z_{\bar{j}}) \right) \\
      & & + i e(\tau^w) \left( i(Z_{\bar{j}})Z_j - i(Z_j) Z_{\bar{j}}
      \right) + i \left( e(\tau^{\bar{j}})Z_{\bar{j}} - e(\tau^j) Z_j
      \right) i(W) \\
      & & + e(\tau^j)i(Z_j)e(\tau^{\bj})i(Z_{\bj})i(W)e(\tau^w) 
      + i(Z_j)e(\tau^j)i(Z_{\bj})e(\tau^{\bj})e(\tau^w)i(W))
\end{eqnarray*}

More simply,
\begin{eqnarray*}
\eta_{j,l} & = & \theta_{j} \theta_{l}^{*} + \theta_{l}^{*} \theta_{j} \\
           & = & \{ -ice(\tau^{j})e(\tau^{\bar{j}})i(W),
                ice(\tau^{w})i(Z_{\bar{l}})i(Z_{l}) \} \\
           & = & c^{2} e(\tau^{j}) e(\tau^{\bar{j}}) i(Z_{\bar{l}})
           i(Z_{l}).
\end{eqnarray*}

      Summing these expressions gives the required formula for the
      Laplacian.

In fact, the equations \eqref{E:ACRfirst}---\eqref{E:ACRlast}, together
with the commutation relations of the Lie group in question, can be
used to define a Lie superalgebra. This theme is developed somewhat in
\cite{thesis}; for more on Lie superalgebras and their connection with
$d$ and the Laplacian, see also \cite{scheunert, sternberg}. 

\section{Proof of the Kernel Lemma}

To prove: if $v \in \ker U_{12} \cap V^{p,n,\gamma}$, then
      $\gamma_2 \le 1$.

      \begin{proof}
      Recall that $U_{12} = a_1^* a_2 - e(\tau^2) i(Z_1) +
      e(\tau^{\bar{1}}) i(Z_{\bar{2}})$.
      Suppose $v \in \ker U_{12}$. Then in particular
      $e(\tau^2) e(\tau^{\bar{1}}) U_{12} v =0$
      which implies that
      $e(\tau^2) e(\tau^{\bar{1}}) a_1^* a_2 v =0$.

      Write $v$ in the form
      $$v= \tau^{\bar{1}} \wedge v_1 + \tau^2 \wedge v_2 +
      \tau^{\bar{1}} \wedge \tau^2 \wedge v_3 + v_4,$$
      for $v_i$ forms such that $i(Z_{\bar{1}})v_i=0=i(Z_2)v_i$
      for $i=1,\dots,4$.
      Then we've just shown above that $v_4 \in \ker a_2$.

      If we apply $U_{12}$ to $v$ and equate coefficients of
      terms with $\tau^2$ and so on, we get the following equations
      (since $v \in \ker U_{12}$):
      \begin{eqnarray}
      a_1^* a_2 v_3 + i(Z_1) v_1 - i(Z_{\bar{2}}) v_2 & = & 0, \label{E:v3} \\
      a_1^* a_2 v_2 - i(Z_1) v_4 & = & 0, \label{E:v2} \\
      a_1^* a_2 v_1 + i(Z_{\bar{2}}) v_4 & = & 0. \label{E:v1}
      \end{eqnarray}

      But we know that $v_4 \in \ker a_2$. Equations \eqref{E:v2} and \eqref{E:v1} then
      imply that $v_2$ and $v_1$ respectively are in $\ker a_2^2$
      (even if $v_4 =0$). From equation \eqref{E:v3}, we see that
      $v_3$ is in $\ker a_2^3$.
      Actually, equation \eqref{E:v1} also implies that $v_1$ is in
      $\ker (i(Z_{\bar{2}}) a_2)$, which together with equation \eqref{E:v3}
      implies that $v_3 \in \ker (i(Z_{\bar{2}}) a_2^2)$.

	If we now require that $v \in V^{p,n,\gamma}$ (and recall that for
	functions, if $\psi_\beta(k) \in \ker a_2^3$, then $\beta_2 \le 2$)
	then the conditions that $\tau^{\bar{1}} \wedge \tau^2 \wedge v_3
	\in V^{p,n,\gamma}$ and $v_3 \in \ker a_2^3 \cap \ker (i(Z_{\bar{2}}) a_2^2)$
	together imply that $\gamma_2 \le 1$, if $v_3 \neq 0$. Similarly the conditions
	on $v_1, v_2$ and $v_4$ imply that $\gamma_2 \le 1$, so that the result holds
	even if one or more of the $v_i$'s is 0.
\end{proof}	

%\bibliography
%\providecommand{\bysame}{\leavevmode\hbox to3em{\hrulefill}\thinspace}


\begin{thebibliography}{10}

\bibitem{atiyah}
M.~Atiyah, \emph{Elliptic operators, discrete groups and von {N}eumann
  algebras}, Ast{\'{e}}risque \textbf{32} (1976), 43--72.

\bibitem{a&s}
M.~Atiyah and W.~Schmid, \emph{A geometric construction of the discrete series
  for semisimple {L}ie groups}, Invent. math. \textbf{42} (1977), 1--62.

\bibitem{bmw}
J.~Block, V.~Mathai, and S.~Weinberger, \emph{Homotopy invariance of
  {N}ovikov-{S}hubin invariants and ${L}^2$ {B}etti numbers}, To appear in
  Proc. Amer. Math. Soc., Nov 1997.

\bibitem{cfm}
A.~Carey, M.~Farber, and V.~Mathai, \emph{Determinant lines, von {N}eumann
  algebras and ${L}^2$ torsion}, J. reine angew. Math. \textbf{484} (1997),
  153--181.

\bibitem{ccmp}
A.L. Carey, T.~Coulhon, V.~Mathai, and J.~Phillips, \emph{Von {N}eumann spectra
  near the spectral gap}, To appear in Bull. Sci. Math. (France).

\bibitem{c&g}
L.J. Corwin and F.P Greenleaf, \emph{Representations of nilpotent {L}ie groups
  and their applications (part 1)}, Cambridge University Press, Cambridge,
  England, 1990.

\bibitem{cdkr}
M.~Cowling, A.~H. Dooley, A.~Kor\'{a}nyi, and F.~Ricci, \emph{${H}$-type groups
  and {I}wasawa decompositions}, Advances in Math. \textbf{87} (1991), 1--41.

\bibitem{dixmierc}
J.~Dixmier, \emph{C*-algebras}, North-Holland, Amsterdam, 1981, Revised
  edition; translation of C*-algebres et leurs representations.

\bibitem{dixmiervn}
\bysame, \emph{Von neumann algebras}, North-Holland, Amsterdam, 1981,
  Translation of Algebres d'operateurs dans l'espace hilbertien (algebres de
  Von Neumann).

\bibitem{dodziuk}
J.~Dodziuk, \emph{De {R}ham-{H}odge theory for ${L}^2$-cohomology of infinite
  coverings}, Topology \textbf{16} (1977), 157--165.

\bibitem{efremov}
A.V. Efremov, \emph{Cellular decompositions and {N}ovikov-{S}hubin invariants},
  Russ. Math. Surveys \textbf{46} (1991), 219--220.

\bibitem{farber}
M.~Farber, \emph{Homological algebra of {Novikov}-{S}hubin invariants and
  {M}orse inequalities}, Geom. Funct. Anal. \textbf{6} (1996), no.~4, 628--665.

\bibitem{folland}
G.B. Folland, \emph{Harmonic analysis in phase space}, Princeton University
  Press, Princeton, N.J., 1989.

\bibitem{follandnew}
\bysame, \emph{A course in abstract harmonic analysis}, CRC Press, Boca Raton,
  c1995.

\bibitem{garding}
L.~G{\aa}rding, \emph{Notes on continuous representations of lie groups}, Proc.
  Nat. Acad. Sci. USA \textbf{33} (1947), 331--332.

\bibitem{g&w}
C.S. Gordon and E.N. Wilson, \emph{The spectrum of the {L}aplacian on
  {R}iemannian {H}eisenberg manifolds}, Mich. Math. J \textbf{33} (1986),
  253--271.

\bibitem{g&s}
M.~Gromov and M.A. Shubin, \emph{Von {N}eumann spectra near zero}, Geom. Anal.
  and Funct. Anal. \textbf{1} (1991), 375--404.

\bibitem{guill&stern}
V.~Guillemin and S.~Sternberg, \emph{Symplectic techniques in physics},
  Cambridge University Press, Cambridge; New York, 1984.

\bibitem{kaplan2}
A.~Kaplan, \emph{On the geometry of {L}ie groups of {H}eisenberg type}, Bull.
  London Math. Soc. \textbf{15} (1983), 35--42.

\bibitem{kirillov}
A.A. Kirillov, \emph{Elements of the theory of representations},
  Springer-Verlag, Berlin; New York, 1976, Translation of Elementy teorii
  predstavlenii.

\bibitem{l&m}
H.B. Lawson, Jr. and M.~Michelsohn, \emph{Spin geometry}, Princeton University
  Press, Princeton, N.J., 1989.

\bibitem{lott}
J.~Lott, \emph{Heat kernels on covering spaces and topological invariants}, J.
  Diff. Geom. \textbf{35} (1992), 471--510.

\bibitem{l&l}
J.~Lott and W.~L{\"{u}}ck, \emph{${L}^2$-topological invariants of
  $3$-manifolds}, Invent. math. \textbf{120} (1995), 15--60.

\bibitem{vmdirac}
V.~Mathai, \emph{Von {N}eumann algebra invariants of {D}irac operators}, To
  appear in J. of Funct. Anal., 1997.

\bibitem{c&m_acyclic}
V.~Mathai and A.~Carey, \emph{${L}^2$-{A}cyclicity and ${L}^2$-torsion
  invariants}, Contemp. Math. \textbf{105} (1990), 91--118.

\bibitem{mautner}
F.I. Mautner, \emph{Unitary representations of locally compact groups}, Annals
  of Math. \textbf{52} (1950), no.~3, 528--556.

\bibitem{ns2}
S.~Novikov and M.A. Shubin, \emph{Morse inequalities and von {N}eumann
  invariants of non-simply-connected manifolds}, Uspekhi Mat. Nauk \textbf{41}
  (1986), no.~5, 222--223, (Russian).

\bibitem{ns1}
\bysame, \emph{Morse theory and von {N}eumann {II}${}_1$-factors}, Doklady
  Akad. Nauk SSSR \textbf{289} (1986), 289--292.

\bibitem{pesce}
H.~Pesce, \emph{Calcul du spectre d'une nilvari{\`{e}}t{\`{e}} de rang deux et
  applications}, Trans. Amer. Math. Soc. \textbf{339} (1993), no.~1, 433--461.

\bibitem{ricci}
F.~Ricci, \emph{Harmonic analysis on generalized {H}eisenberg groups},
  Preprint.

\bibitem{satake}
I.~Satake, \emph{Linear algebra}, M. Dekker, New York, 1975, Translation of
  Senkei daisugaku.

\bibitem{scheunert}
M.~Scheunert, \emph{The theory of {L}ie superalgebras: an introduction},
  Springer-Verlag, Berlin; New York, 1979.

\bibitem{thesis}
L.~Schubert, \emph{Spectral properties of the {L}aplacian on $p$-forms on the
  {H}eisenberg group}, Ph.D. thesis, The University of Adelaide, 1997.

\bibitem{segal}
I.M. Segal, \emph{A non-commutative extension of abstract integration}, Ann.
  Math. \textbf{57} (1953), 401--457, Corrections in 58 (1954), pp595--596.

\bibitem{spivak}
M.~Spivak, \emph{A comprehensive introduction to differential geometry (volume
  1)}, Publish or Perish Inc., Berkeley, 1979.

\bibitem{sternberg}
S.~Sternberg, \emph{Some recent results on the metaplectic representation},
  Group theoretical methods in physics (New York) (P.~Kramer and A.~Rieckers,
  eds.), Lecture Notes in Physics, no.~79, Springer-Verlag, 1978.

\bibitem{varopoulos}
N.~Varopoulos, \emph{Random walks and {B}rownian motion on manifolds}, Sympos.
  Math. \textbf{29} (1988), 97--109.

\end{thebibliography}
\end{document}